\documentstyle{article}

\newtheorem{theorem}{Theorem}[section]

\newtheorem{lemma}[theorem]{Lemma}

\newtheorem{definition}[theorem]{Definition}
\newtheorem{example}[theorem]{Example}

\newcommand{\defin}[1]
  {\begin{definition} {\rm #1} \end{definition}}
\newcommand{\examp}[1]
  {\begin{example} {\rm #1} \end{example}}

\def\QED{\quad\blackslug\lower 8.5pt\null}

\begin{document}

{\Large \bf  Semiintegrable almost Grassmann 
structures}

\vspace*{3mm}

{\large M. A. Akivis\footnote{
 The research of the first author was 
 partially supported by the Israel 
Ministry of Absorption and the Israel Public Council for Soviet 
Jewry.}}  

{\footnotesize \it             
Department of Mathematics,   Jerusalem College of Technology 
- Mahon Lev, 21 Havaad  

Haleumi St.,  
P. O. B. 16031, Jerusalem 91160, Israel;  E-mail address: akivis@math.jct.ac.il}

\vspace*{3mm}

{\large   V. V. Goldberg\footnote{The research of the second 
author  was partially 
supported by the Research Council of the Catholic University of Leuven, Leuven, Belgium.}} 

{\footnotesize \it             
Department of Mathematics,   New Jersey Institute of Technology,   University Heights, 

Newark, NJ 07102, U.S.A.;  E-mail address: vlgold@numerics.njit.edu}

\vspace*{2mm}

{\footnotesize 
{\it Abstract:} 
In the present paper we study locally semiflat (we also call 
them semiintegrable) 

almost Grassmann structures. We establish necessary and 
sufficient conditions for an 

almost Grassmann structure to be $\alpha$- or 
$\beta$-semiintegrable.  These conditions are 

expressed in terms of the fundamental tensors of almost Grassmann 
structures. Since we 
are not able to prove the existence of locally semiflat 
almost Grassmann structures in 
 the general case, we give many examples of $\alpha$- and 
$\beta$-semiintegrable structures, mostly 
four-dimensional. For all examples we find systems 
of differential equations of the 
 families of integral submanifolds $V_\alpha$ and $V_\beta$ of 
the distributions  $\Delta_\alpha$ and $\Delta_\beta$ of  plane 
elements associated with an almost Grassmann 
structure. For some examples we  were 
able to integrate these systems and find closed 
form equations of submanifolds $V_\alpha$ 
and $V_\beta$.

\vspace*{2mm}


\vspace*{2mm}

{\it MS classification}: 53A40, 53A60.
}

\setcounter{equation}{0}

\setcounter{section}{-1}

\section{Introduction}  In the paper 
\cite{AG97} (see also the book \cite{AG96}, Ch. 7) the authors 
constructed the real theory of almost Grassmann structures 
$AG (p-1, p+q-1)$  defined 
on a differentiable manifold of dimension $n = pq$ by a fibration 
of Segre cones $SC (p, q)$.  In particular, in \cite{AG97} we  
derived the structure equations of $AG (p-1, p+q-1)$ and found (in a fourth-order differential neighborhood)   a complete  
geometric 
object of the almost Grassmann  structure totally defining its 
geometric structure. We also  found the structure group of   
these structures and its  
differential prolongation and  the conditions 
under which an almost Grassmann structure is locally flat or 
locally semiflat.

While constructing this theory, we distinguished three 
cases:  $p = 2, \; q = 2 \; (\dim M = 4)$;  $p = 2, \; q > 2$ (or 
$p > 2, \; q = 2$); and  $p > 2, \; q > 2$. 
 We constructed the fundamental geometric objects of these 
structures up to fourth order for each of these three cases and 
established connections among them.
In the first  case the almost Grassmann  structure 
$AG (1, 3)$ is equivalent to the pseudoconformal structure 
$CO (2, 2)$.  Since the 
four-dimensional conformal structures play an important 
role in general relativity, this provides a physical 
justification for   studying  the structures $AG (1, 3)$ 
as well as for studying  
the general almost Grassmann structures $AG (p-1, p+q-1)$.

In the present paper we study locally semiflat almost Grassmann 
structures (we also call them  semiintegrable), 
and for different values $p$ and $q$ we establish necessary and 
sufficient conditions of $\alpha$- and $\beta$-semiintegrability 
of almost Grassmann  structures. These conditions are expressed 
in terms of the fundamental tensors of almost Grassmann structures. 
We find the relations between  10 independent  components 
of  an almost Grassmann  structure $AG (1, 3)$ and 
10 independent  components of an equivalent pseudoconformal structure 
$CO (2, 2)$.

Since we are not able to prove the existence of locally semiflat 
almost Grassmann structures in the general case, we 
give many examples of $\alpha$- and $\beta$-semiintegrable 
structures, mostly for $p = q = 2$. For all examples we find 
systems of differential equations of 
the families of integral submanifolds $V_\alpha$ and $V_\beta$ 
of the distributions $\Delta_\alpha$ and $\Delta_\beta$ 
of plane elements 
associated with almost Grassmann structures. For some examples we 
were able to integrate these systems and find closed form 
equations of submanifolds $V_\alpha$ and $V_\beta$. 

Note that the existence of globally semiflat four-dimensional 
conformal structures was proved in \cite{T92} (see also 
\cite{LB95}).

\section{Almost Grassmann Structures}

\setcounter{equation}{0}

{\bf 1.} 
First we recall the definition of Segre varieties and 
 Segre cones.

The  {\em Segre variety} $S (k, l)$ is an embedding of the direct 
product $P^k \times P^l$ of projective spaces $P^k$ and $P^l$ 
of dimensions $k$ and $l$ into a projective space of dimension 
$(k+1)(l+1)-1 = kl + k + l$. Analytically this embedding can be 
written by means of the following equations:
\begin{equation}\label{eq:1.1}
z_\alpha^i = t_\alpha s^i, \;\;\;\;\; 
\alpha = 0, 1, \ldots, k; \;\; i = 0, 1, \ldots, l,
\end{equation}
where $t_\alpha, s^i$ and $z_\alpha^i$ are homogeneous 
coordinates in the spaces 
$P^k, P^l$ and $P^{kl+k+l}$, respectively. These equations are equivalent to 
the condition
\begin{equation}\label{eq:1.2}
\mbox{rank} \, (z_\alpha^i) = 1.
\end{equation}
The Segre variety $S (k, l)$ has the dimension $k+l$. 
It is proved in algebraic geometry that the degree of this 
variety is 
$$
\mbox{deg} \;\; S(k, l) = {k+l \choose k}.
$$

The cone $SC_x (k+1, l+1)$ with vertex at the point $x$ 
whose projectivization is the Segre 
variety  $S (k, l)$ is called the {\em Segre cone}. 

Now we can define the notion of almost Grassmann structure. 

\defin{\label{def:1.1} Let $M$ be a  differentiable manifold of 
dimension $p q$, and let $SC (p, q)$ be a differentiable 
fibration of Segre cones with the base $M$ such that 
$SC_x (M) \subset T_x (M), \;\; x \in M$. The pair 
$(M, SC (p, q))$ is said to be  an {\it almost Grassmann 
structure} and is denoted by $AG (p-1, p+q-1)$. The manifold 
$M$ endowed with  such a structure is said to be an {\it almost 
Grassmann manifold}.
}

Note that the almost Grassmann 
structure $AG (p-1, p+q-1)$ is equivalent to the structure 
$AG (q-1, p+q-1)$ since both of these structures are generated on 
the manifold $M$ by a differentiable family of  Segre cones 
$SC_x (p, q)$. 

In \cite{AG97} we discussed the following 
examples of  almost Grassmann structures: the almost 
Grassmann structure associated with the Grassmannian $G (m, n)$ 
(in this case $p = m+1$ and $q = n-m$); the almost 
Grassmann structure $AG (1, 3)$ which is equivalent to 
the pseudoconformal $CO (2, 2)$-structure; and almost Grassmann 
associated with multidimensional webs.

{\bf 2.} 
The structural group of the almost Grassmann structure is a 
subgroup of the general linear group ${\bf GL}(pq)$ of 
transformations of the space $T_{x} (M)$, which leave the cone 
$SC_{x}(p, q) \subset T_x (M)$ invariant. We denote this group by 
$G = {\bf GL} (p, q)$. 

To clarify the structure of this group,  we consider in the tangent space  $T_x (M)$ a family of frames  
$\{e_{i}^{\alpha}\}, \alpha = 1, \ldots, p; \;\; i=p+1, \ldots, 
p+q$, such that for any fixed $i$, the vectors $e_i^\alpha$ 
belong to a $p$-dimensional generator $\xi$ of the Segre cone 
$SC_x (p, q)$, and for any fixed $\alpha$, the vectors 
$e_i^\alpha$ belong to a $q$-dimensional generator 
$\eta$ of  $SC_x (p, q)$.  In such a frame, the equations  
of the cone $SC_x (p, q)$ can be written in the form 
(1.1) where now $\alpha = 1, \ldots, p; \; i = p+1, \ldots , p+q, 
z_\alpha^i$ are the coordinates of a vector 
$z = z^i_\alpha e_i^\alpha \subset T_x (M)$, and $t_\alpha$ and 
$s^i$ are parameters on which a vector $z \subset SC_x (M)$ 
depends. 

As was shown in \cite{AG96}, the group $G$ of admissible 
transformations of the frames $\{e_i^\alpha\}$ keeping the Segre 
cone $SC_x (p, q)$ invariant can be presented in the form:
\begin{equation}\label{eq:1.3}
G = {\bf SL} (p) \times {\bf SL}(q)  \times {\bf H},
\end{equation}
where ${\bf SL} (p)$ and ${\bf SL} (q)$ are special linear groups 
in spaces of dimensions $p$ and $q$, and 
${\bf H} = {\bf R}^* \otimes \mbox{{\rm Id}}$ is the group  
of homotheties in  $T_x (M)$. 

It follows that {\em an almost Grassmann structure  $AG (m, n)$ 
is a $G$-structure of first order.}

 From equation (1.1)  defining the Segre cone $SC_x (p, q)$  it follows that 
this cone carries $(q - 1)$-parameter family of 
$p$-dimensional generators $\xi$  and $(p-1)$-parameter family 
of $q$-dimensional generators $\eta$.

The $p$-dimensional generators $\xi$ form a fiber bundle on 
the manifold $M$. The base of this bundle is the manifold $M$, 
and its fiber attached to a point $x \in M$ is the set of all 
$p$-dimensional plane generators $\xi$ of the Segre cone 
$SC_x (p, q)$. The dimension of a fiber is $q - 1$, and it is 
parametrized by means of a projective space 
$P_\alpha, \; \dim P_\alpha = q - 1$. We will denote this fiber 
bundle of $p$-subspaces by $E_\alpha = (M, P_\alpha)$. 

In a similar manner, $q$-dimensional plane generators $\eta$ of 
the Segre cone $SC_x (p, q)$ form on $M$ the fiber bundle 
$E_\beta = (M, P_\beta)$ with the base $M$ and fibers of 
dimension $p - 1 = \dim P_\beta$. The fibers are $q$-dimensional 
plane generators $\eta$ of the Segre cone $SC_x (p, q)$. 

Consider the manifold $M_\alpha = M \times P_\alpha$ of dimension 
$pq + q - 1$. The fiber bundle $E_\alpha$ induces on $M_\alpha$ 
the distribution $\Delta_\alpha$ of plane elements $\xi_\alpha$ 
of dimension $q$. In a similar manner, on the manifold 
$M_\beta = M \times P_\beta$ the fiber bundle $E_\beta$ induces 
the distribution $\Delta_\beta$ of plane elements $\xi_\beta$ of 
dimension $p$. 

\defin{\label{def:1.2}
An almost Grassmann structure $AG (p-1, p+q-1)$ is said to be 
{\em $\alpha$-semiintegrable} if the distribution $\Delta_\alpha$ 
is  integrable on this structure.  Similarly, an almost Grassmann 
structure $AG (p-1, p+q-1)$ is said to be 
{\em $\beta$-semiintegrable} if  the distribution $\Delta_\beta$ 
is  integrable on this structure. A structure $AG (p-1, p+q-1)$ 
is called {\em integrable} if it is both $\alpha$- and 
$\beta$-semiintegrable.
} 

Integral manifolds $\widetilde{V}_\alpha$ of the distribution 
$\Delta_\alpha$ of an $\alpha$-semiintegrable almost Grassmann 
structure are of dimension $p$. They are projected onto the 
original manifold $M$ in the form of a submanifold $V_\alpha$ of 
the same dimension $p$, which, at any of its points, is tangent 
to the $p$-subspace $\xi_\alpha $ of the fiber bundle $E_\alpha$. 
Through each point $x \in M$, there passes a $(q-1)$-parameter 
family of submanifolds $V_\alpha$.

Similarly, integral manifolds $\widetilde{V}_\beta$ of the 
distribution $\Delta_\beta$ of a $\beta$-semiinte\-gra\-ble 
almost 
Grassmann structure are of dimension $q$. They are projected onto 
the original manifold $M$ in the form of a submanifold $V_\beta$ 
of the same dimension $q$, which, at any of its points, is 
tangent to the $q$-subspace $\eta_\beta$ of the fiber bundle 
$E_\beta$. Through each point $x \in M$, there passes a 
$(p-1)$-parameter family of submanifolds $V_\beta$. 
If an almost Grassmann structure on $M$ is integrable, then 
through each point $x \in M$, there pass a $(q-1)$-parameter 
family of submanifolds $V_\alpha$ and a $(p-1)$-parameter family 
of submanifolds $V_\beta$ which were described above. 

The Grassmann structure $G (m, n)$ is an integrable almost 
Grassmann structure $AG (m, n)$ since through 
any point $x$ of the variety $ \Omega (m, n)$, onto which the manifold 
$G (m, n)$ is mapped bijectively  under the Grassmann mapping, 
there pass a $(q - 1)$-parameter family of $p$-dimensional plane 
generators (which are the submanifolds $V_\alpha$) and a 
$(p - 1)$-parameter family of $q$-dimensional plane generators 
(which are the submanifolds $V_\beta$). In the projective space 
$P^n$,  to  submanifolds $V_\alpha$ there corresponds a family of 
$m$-dimensional subspaces belonging to a subspace of dimension 
$m + 1$, and   to  submanifolds $V_\beta$ there corresponds a 
family of $m$-dimensional subspaces passing through a subspace of 
dimension $m - 1$.

Note that if the manifold $M$ orientable, and  
we change its orientation, then an 
$\alpha$-semiintegrable almost Grassmann structure will 
become $\beta$-integrable, and vice versa. 

{\bf 3.} 
Consider  a differentiable manifold $M$ of dimension $pq$ 
endowed with an almost Grassmann structure 
$AG (p-1, p+q-1)$. Suppose that $x \in M$,  $T_x (M)$ 
is the tangent space of the manifold    $M$ at the 
point $x$ and that $\{e_i^\alpha\}$ is an adapted frame 
of the structure $AG (p-1, p+q-1)$. The decomposition of 
a vector $z \in T_x (M)$ with respect to this basis can be 
written in the form
$$
z = \omega_\alpha^i (z) e_i^\alpha,
$$
where $\omega_\alpha^i$ are 1-forms making up the {\em co-frame} 
in the space $T_x (M)$. If $z = dx$ is the differential of a point $x \in M$, 
then the forms $\omega_\alpha^i (dx)$ are 
differential forms defined on a first-order frame bundle 
associated with the almost Grassmann structure. These forms 
constitute a completely integrable system of forms. 

As was proved in \cite{AG97}, the form 
$\theta = (\omega_\alpha^i)$ 
and the forms arising under its prolongation 
satisfy the following structure equations:
\begin{equation}\label{eq:1.4}
\renewcommand{\arraystretch}{1.5}
\left\{
 \begin{array}{ll}
 d\omega_\alpha^i - \omega_\alpha^j \wedge 
\omega_j^i - \omega_\alpha^\beta \wedge \omega_\beta^i 
-  \omega \wedge \omega_\alpha^i = a_{\alpha jk}^{i\beta\gamma} 
\omega_{\beta}^{j} \wedge \omega_{\gamma}^{k}, \\
d \omega_\alpha^\beta - \omega_\alpha^\gamma \wedge \omega_\gamma^\beta 
= \displaystyle \frac{q}{p+q}\Bigl(\delta_\alpha^\beta 
\omega_\gamma^k \wedge \omega^\gamma_k - p \omega_\alpha^k 
\wedge \omega_k^\beta\Bigr) 
+  b_{\alpha kl}^{\beta\gamma\delta}  \omega_\gamma^k \wedge 
\omega_\delta^l, \\
 d \omega_j^i - \omega_j^k  \wedge \omega_k^i
= \displaystyle \frac{p}{p+q} \Bigl(\delta_j^i \omega^\gamma_k 
\wedge \omega_\gamma^k 
- q \omega_j^\gamma  \wedge \omega_\gamma^i\Bigr) 
+  b_{jkl}^{i\gamma\delta}  \omega_\gamma^k \wedge 
\omega_\delta^l, \\
d \omega = \omega_i^\alpha \wedge \omega_\alpha^i,\\
d \omega_i^\alpha - \omega_i^\beta \wedge \omega_\beta^\alpha 
- \omega_i^j \wedge \omega_j^\alpha  
+ \omega \wedge \omega_i^\alpha
=    c_{ijk}^{\alpha\beta\gamma} 
 \omega^k_\gamma \wedge \omega_\beta^j 
- a_{\gamma ij}^{k\alpha\beta} \omega_k^\gamma \wedge 
\omega_\beta^j,
 \end{array}
\right.
\renewcommand{\arraystretch}{1}
\end{equation}
where  
the matrix 1-form $\theta = (\omega_\alpha^i)$ is defined in 
a first-order frame bundle,  the  form $\omega$ 
is a scalar form   defined in a second-order frame bundle,  
$\omega_\alpha = (\omega_\beta^\alpha)$ and 
$\omega_\beta = (\omega_j^i)$ are the matrix 1-forms also 
defined in a second-order frame bundle. 

The forms $\omega_\alpha^\beta$ and $\omega_i^j$ 
satisfy the conditions
\begin{equation}\label{eq:1.5}
\omega_\gamma^\gamma = 0, \;\; \omega_k^k = 0.
\end{equation}

In equations (1.4) the 2-form 
\begin{equation}\label{eq:1.6}
\Theta^i_\alpha = a_{\alpha jk}^{i\beta\gamma} 
\omega_{\beta}^{j} \wedge \omega_{\gamma}^{k}
\end{equation}
is the {\em torsion $2$-form} of the $AG (p-1, p+q-1)$-structure, 
and the forms 
\begin{equation}\label{eq:1.7}
 \Omega_\alpha^\beta
=  b_{\alpha kl}^{\beta\gamma\delta}  \omega_\gamma^k \wedge 
\omega_\delta^l, \;\;
 \Omega_j^i 
=  b_{jkl}^{i\gamma\delta}  \omega_\gamma^k \wedge 
\omega_\delta^l, \;\;
 \Phi_i^\alpha =   c_{ijk}^{\alpha\beta\gamma} 
 \omega^k_\gamma \wedge \omega_\beta^j 
\end{equation}
are the {\em curvature $2$-forms} of this structure. 

Note that Dhooghe in  \cite{D93} and \cite{D94} 
 gave the structure equations of an almost Grassmann 
structure in the form close to our equations (1.4). His equations 
differ from equations (1.4) only by the additional 
term $\Omega_*$ in the right-hand 
side of the fourth equation. This means that on a manifold 
$M, \; \dim M = pq$, carrying 
an almost Grassmann structure our equations (1.4) 
define a {\em normal Cartan connection}. 

Moreover, in his further considerations in 
\cite{D93} and \cite{D94} Dhooghe assumes that the torsion form 
$\Omega_\alpha^i$ is identically equal to 0 not 
only for $ p = q = 2$ but also for $p > 2$ and $q > 2$. 
This leads to the loss of generality. We did not make 
the above assumption.

In \cite{AG97} (see also \cite{AG96}, 
\S 7.2) we proved the following facts: 

\begin{description}
\item[a)] 
The quantities $ a = \{a^{i\beta\gamma}_{\alpha jk}\}$, 
defined by a second-order neighborhood, 
 form a relative tensor of weight $- 1$ 
and satisfy the following conditions:
\begin{equation}\label{eq:1.8}
a^{i\beta\gamma}_{\alpha jk} = - a^{i\gamma\beta}_{\alpha kj}
\end{equation}
and
\begin{equation}\label{eq:1.9}
a^{i\alpha\gamma}_{\alpha jk} = 0, \;\;\;
a^{i\beta\gamma}_{\alpha ik} = 0.
\end{equation}
The tensor $\{a_{\alpha jk}^{i\beta\gamma}\}$ is said to be the 
{\it  first structure tensor}, or the {\it torsion tensor}, of an 
almost Grassmann manifold $AG (p-1, p+q-1)$.

\item[b)] 
Let us set $b^1 = \{b_{jkl}^{i\gamma\delta}\}, 
b^2 = \{b_{\alpha kl}^{\beta\gamma\delta}\}$ and 
$b = (b^1, b^2)$. The quantities 
$(a, b^1)$ and $(a, b^2)$ form linear homogeneous 
objects. They represent two subobjects of the {\em second 
structure object} $(a, b)$  of the  almost Grassmann 
structure $AG (p-1, p+q-1)$. The components 
$b^1$ and $b^2$  satisfy the conditions:
\begin{equation}\label{eq:1.10}
  b_{\alpha kl}^{\beta\gamma\delta} 
= - b_{\alpha lk}^{\beta\delta\gamma}, \;\; 
  b_{jkl}^{i\gamma\delta} = - b_{jlk}^{i\delta\gamma},
\end{equation}
\begin{equation}\label{eq:1.11}
  b_{\alpha kl}^{\alpha\gamma\delta} = 0,  \;\; 
  b_{ikl}^{i\gamma\delta} = 0,
\end{equation}
\begin{equation}\label{eq:1.12}
 b^{\gamma\alpha\delta}_{\alpha kl} 
-  b^{i\gamma\delta}_{kil}  
+ b^{\delta\alpha\gamma}_{\alpha lk} 
- b^{i\delta\gamma}_{lik} = 0,
\end{equation}
and the components of the tensor $a$ satisfy the 
following differential equations:
\begin{equation}\label{eq:1.13}
2 a_{\alpha [j|m}^{i [\beta|\varepsilon} 
a_{\epsilon|kl]}^{m|\gamma\delta]} + \delta_{[j}^{i} 
b_{|\alpha|kl]}^{[\beta\gamma\delta]} 
- \delta_{\alpha}^{[\beta} b_{[jkl]}^{|i|\gamma\delta]} 
+ a_{\alpha [jkl]}^{i[\beta\gamma\delta]} = 0, 
 \end{equation}
where $a^{i\beta\gamma\delta}_{\alpha jkl}$  
are the Pfaffian derivatives of $a^{i\beta\gamma}_{\alpha jk}$ 
and the alternation is performed with respect to 
the vertical pairs of indices ${\beta \choose j}$ and ${\gamma
 \choose k}$. Formulas (1.13) are analogues of the Bianchi 
equations in the theory of spaces with affine connection.

\item[c)]  For $p > 2$ and $q > 2$, 
 the components of $b^2$ and $b^1$, respectively, are 
expressed in terms of the components of the tensor $a$ and their 
Pfaffian derivatives. 

This implies that if for $p > 2, q > 2$, 
  the torsion tensor $a$  of an almost Grassmann structure 
vanishes, then its curvature tensor $b$ vanishes as well.

\item[d)] 
Let us set $c = \{c^{\alpha\beta\gamma}_{ijk}\}$. 
Then  $S = (a, b, c)$ 
 forms a linear homogeneous object, which is called the {\em 
third structure object} of the almost Grassmann structure 
$AG (p-1, p+q-1)$. It is defined by a fourth-order differential 
neighborhood of $AG (p-1, p+q-1)$. 
Its subobject $a$ is a relative tensor (the torsion tensor) 
 defined by a second-order differential 
neighborhood of $AG (p-1, p+q-1)$, and the subobjects 
$(a, b^1), (a, b^2)$, and $(a, b)$ are defined by a third-order 
differential neighborhood of $AG (p-1, p+q-1)$. 

 The third structural object $S = (a, b, c)$ 
is the {\em complete geometric object} of the almost Grassmann 
structure $AG (p - 1, p + q - 1)$, since during the prolongation  
of  structure equations (1.4) of $AG (p - 1, p + q - 1)$, 
all newly arising objects are expressed in terms of 
the components of the object $S$ and their Pfaffian derivatives.  

The components of  $c$  satisfy the conditions 
\begin{equation}\label{eq:1.14}
c^{i\beta\gamma}_{\alpha jk} = - c^{i\gamma\beta}_{\alpha kj}
\end{equation}
and
\begin{equation}\label{eq:1.15}
  c_{[ijk]}^{[\alpha\beta\gamma]} = 0,
\end{equation}
and the components of $b$ satisfy the  
differential equations
\begin{equation}\label{eq:1.16}
\renewcommand{\arraystretch}{1.5}
\left\{
 \begin{array}{ll}
 b^{\beta[\gamma\delta\varepsilon]}_{\alpha [klm]} 
- \displaystyle \frac{pq}{p+q}  \delta_\alpha^{[\varepsilon}
  c_{[mkl]}^{|\beta|\gamma\delta]}
- 2 b^{\beta[\gamma|\sigma}_{\alpha [k|s} 
a^{s| \delta\varepsilon]}_{\sigma| lm]} = 0,\\
 b^{i[\gamma\delta\varepsilon]}_{j[klm]} 
+ \displaystyle \frac{pq}{p+q}  
  \delta_{[m}^i c^{[\varepsilon\gamma\delta]}_{|j|kl]} 
+ 2 b^{i\sigma[\varepsilon}_{js[m} 
a^{|s| \gamma\delta]}_{|\sigma| kl]} = 0,
\end{array} 
\right.
\renewcommand{\arraystretch}{1}
\end{equation}
where $b^{i\gamma\delta\varepsilon}_{jklm}$ and 
$b^{\beta\gamma\delta\varepsilon}_{\alpha klm}$ 
are the 
Pfaffian derivatives of $b^{\beta\gamma\delta}_{\alpha kl}$ 
and $b^{\beta\gamma\delta}_{\alpha kl}$, respectively.
In formulas (1.15) and (1.16) the alternation is carried out 
with respect to the vertical pairs of indices.

\item[e)] 
If $p > 2$, then the components of $c$  are 
expressed in terms of the components of the subobject $(a, b^1)$ 
and their Pfaffian derivatives, and  if $q > 2$, then 
the components of $c$  are expressed in terms of the components 
of the subobject $(a, b^2)$ and their Pfaffian derivatives. 
This implies that in this case  the object $(a, b)$  
satisfies certain differential equations which 
are other analogues of the Bianchi 
equations in the theory of spaces with affine connection.
 These equations can be obtained if we substitute for the 
components  of $c$  in equations (1.16)  their values.

\item[f)] 
An almost  Grassmann structure $AG (p-1, p+q-1)$ 
is said to be {\em locally Grassmann} (or {\em locally flat}) 
if it is locally equivalent to a Grassmann structure. 
 {\em For $p > 2$ and $q > 2$,  an almost  Grassmann structure  
 $AG (p-1, p+q-1)$ is locally Grassmann if and only if 
its first structure tensor $a$ vanishes. For $p = 2, q = 2$, the tensor 
$a^{i\beta\gamma}_{\alpha jk}$  vanishes, and the condition 
for the structure $AG (1, 3)$ to be locally Grassmann is 
the vanishing of the tensor $b$}. The case $p = 2, q>2$ 
(and $p>2, q=2$) will be considered below.
\end{description}

\setcounter{equation}{0}

\section{Semiintegrability of Almost Grassmann Structures}

{\bf 1.} In this section we  find  geometric conditions for 
an almost Grassmann structure $AG(p-1, p+q-1)$ defined on a 
manifold $M$ to be semiintegrable. The conditions are 
expressed in terms of the complete structure object 
$S$ of  the almost Grassmann structure 
$AG (p-1, p+q-1)$ and its subobjects $S_\alpha$ and $S_\beta$ 
which will be defined in this section. 

In what follows, we 
 often encounter quantities satisfying the conditions 
similar to conditions (1.8) for the quantities 
$a^{i\beta\gamma}_{\alpha jk}$. For calculations with 
quantities of this kind, the following lemma   
is very useful: 

\begin{lemma} 
If a system of quantities $T^{\ldots \alpha \beta}_{\ldots ij}$ 
is skew-symmetric with respect to 
the pairs of indices $
{\alpha \choose i}$ and
 $
{\beta \choose j}$, namely satisfies 
the conditions 
\begin{equation}\label{eq:2.1}
T^{\ldots \alpha \beta}_{\ldots ij} 
= - T^{\ldots \beta \alpha }_{\ldots ji}, 
\end{equation}
then the following identities hold:
\begin{equation}\label{eq:2.2}
\renewcommand{\arraystretch}{1.5}
\begin{array}{ll}
T^{\ldots [\alpha \beta]}_{\ldots ij} 
= T^{\ldots \alpha \beta}_{\ldots (ij)}, &
T^{\ldots (\alpha \beta)}_{\ldots ij} 
= T^{\ldots \alpha \beta}_{\ldots [ij]}, \\
T^{\ldots (\alpha \beta)}_{\ldots ij} 
= - T^{\ldots (\alpha \beta)}_{\ldots ji} 
= T^{\ldots (\alpha \beta)}_{\ldots [ij]}, &
T^{\ldots \alpha \beta}_{\ldots [ij]} 
= - T^{\ldots  \beta \alpha}_{\ldots [ij]} 
= T^{\ldots (\alpha \beta)}_{\ldots [ij]}, \\
T^{\ldots [\alpha \beta]}_{\ldots [ij]} = 0, & 
T^{\ldots (\alpha \beta)}_{\ldots (ij)} = 0.
\end{array}
\renewcommand{\arraystretch}{1}
\end{equation}
In these relations the symmetrization and the alternation
 are carried separately over  the lower 
indices and the upper indices. In addition the following 
decompositions take place:
\begin{equation}\label{eq:2.3}
T^{\ldots \alpha \beta}_{\ldots ij} 
= T^{\ldots (\alpha \beta)}_{\ldots ij}
+ T^{\ldots \alpha \beta}_{\ldots (ij)}, \;\; 
 T^{\ldots \alpha \beta}_{\ldots ij}
= T^{\ldots [\alpha \beta]}_{\ldots ij} 
+  T^{\ldots \alpha \beta}_{\ldots [ij]}. 
\end{equation}
\end{lemma}

{\sf Proof.} All these identities can be proved by 
direct calculation with help of (2.1). \rule{3mm}{3mm}

In addition, in the proof of the main theorem, we will use the 
following lemma: 

\begin{lemma} 
The condition 
\begin{equation}\label{eq:2.4}
T^{[\alpha\beta\gamma]}_{[ijk]} = 0,
\end{equation} 
where the alternation is carried over three vertical pairs of 
indices,  implies the condition 
\begin{equation}\label{eq:2.5}
T^{[\alpha\beta\gamma]}_{(ijk)} = 0,
\end{equation} 
where the alternation and symmetrization are carried  
separately over the upper triple of indices and the lower triple 
of indices.
\end{lemma}

{\sf Proof.} To prove this, one writes down 
36 terms of  $T^{[\alpha\beta\gamma]}_{(ijk)}$ and 
collects from them 6 groups of 6 terms to each of which 
the hypothesis (2.4) can be applied. \rule{3mm}{3mm}

Next we will prove the following important result 
on the decomposition of the torsion tensor of an almost Grassmann 
structure $AG (p-1, p+q-1)$: 

\begin{theorem}
The torsion tensor $a  \! 
= \! \{a_{\alpha jk}^{i\beta\gamma}\}$ 
of the almost Grassmann structure 
$AG (p-1, p+q-1)$  decomposes into two subtensors:
\begin{equation}\label{eq:2.6}
a \!= \!a_\alpha \dot{+} a_\beta,
\end{equation}
where 
$$
a_\alpha =  \{a_{\alpha (jk)}^{i\beta\gamma}\}, \;\; 
a_\beta = \{ a_{\alpha jk}^{i(\beta\gamma)}\}.
$$
\end{theorem}

{\sf Proof.} Since the tensor 
$a_{\alpha jk}^{i\beta\gamma}$ is skew-symmetric with respect to 
the pairs of indices  $
{\beta \choose j}$ and 
$
{\gamma \choose k}$, then,  by Lemma 2.1, 
the decomposition (2.6) is equivalent to the obvious 
decomposition 
$$
 a_{\alpha jk}^{i\beta\gamma} 
= a_{\alpha (jk)}^{i\beta\gamma} 
+ a_{\alpha [jk]}^{i\beta\gamma}.  \rule{3mm}{3mm}
$$

Note that by Lemma 2.1, the subtensors $a_\alpha$ and $a_\beta$ 
can be also represented in the form 
$$
a_\alpha =  \{a_{\alpha jk}^{i[\beta\gamma]}\}, \;\; 
a_\beta = \{ a_{\alpha [jk]}^{i\beta\gamma}\}.
$$

Note also that like the tensor $a$, its subtensors $a_\alpha$ and 
$a_\beta$ are skew-symmetric  with respect to 
the pairs of indices  $
{\beta \choose j}$ and 
$
{\gamma \choose k}$:
$$
a_{\alpha (jk)}^{i\beta\gamma} = -  
a_{\alpha (kj)}^{i\gamma\beta}, \;\; 
a_{\alpha [jk]}^{i\beta\gamma} = -  
a_{\alpha [kj]}^{i\gamma\beta},  
$$

\noindent
  and they are also trace-free, 
since it follows from  (1.9) that 
\begin{equation}\label{eq:2.7}
a_{\alpha (jk)}^{i\alpha\gamma} = 0, \;\;
a_{\alpha ik}^{i[\beta\gamma]} = 0, \;\; 
a_{\alpha ik}^{i(\beta\gamma)} = 0, \;\;
a_{\alpha [jk]}^{i\alpha\gamma} = 0. 
\end{equation}

\begin{theorem}
If $p = 2$, then $a_\alpha  = 0$, 
and if $q = 2$, then $a_\beta = 0$. 
\end{theorem}

{\sf Proof.} Suppose that $p = 2$. Then $\alpha, \beta, \gamma 
= 1, 2$. Since, by definition and Lemma 2.1, the tensor 
$a_\alpha$ is skew-symmetric with respect to the indices 
$\beta$ and $\gamma$, we have
$$
a_{\alpha (jk)}^{i 11} =  a_{\alpha (jk)}^{i 22} = 0.
$$
But the first condition of (2.7) gives
$$
a_{1 (jk)}^{i 11} + a_{2 (jk)}^{i 21} = 0,\;\; 
a_{1 (jk)}^{i 12} + a_{2 (jk)}^{i 22} = 0. 
$$
It follows from these relations that 
 $$
a_{2 (jk)}^{i 21} =  a_{1 (jk)}^{i 12} = 0;
$$
 that is, all components of the tensor $a_\alpha$ vanish. 

For the case $q=2$, the proof is similar.  \rule{3mm}{3mm}
\vspace*{2mm}

{\bf 2.} We introduce the following notation:
\begin{equation}\label{eq:2.8}
\renewcommand{\arraystretch}{1.3}
\begin{array}{lll}
b_\alpha^1 = \{b_{(jkl)}^{i\gamma\delta}\}, & \!\! 
b_\alpha^2 = \{b_{\alpha kl}^{[\beta\gamma\delta]}\}, & \!\!
  c_\alpha = \{c_{(ijk)}^{[\alpha\beta\gamma]}\}, \\ 
&&\\
b_\beta^1 = \{b_{[jkl]}^{i\gamma\delta}\}, & \!\! 
b_\beta^2 = \{b_{\alpha kl}^{(\beta\gamma\delta)}\}, & \!\! 
  c_\beta = \{c_{[ijk]}^{(\alpha\beta\gamma)}\}. 
\end{array}
\renewcommand{\arraystretch}{1}
\end{equation}

Now we give necessary and sufficient conditions for an 
almost Grassmann structure $AG (p-1, p+q-1)$ to 
be $\alpha$-semiintegrable or  $\beta$-semiintegrable.

\begin{theorem}

\begin{description}
\item[(i)] 
If $p > 2$ and $q \geq 2$, then for an almost Grassmann structure 
$AG (p-1, p+q-1)$ to be $\alpha$-semiintegrable, it is necessary 
and sufficient that the conditions    $a_\alpha = b_\alpha^1 
= b_\alpha^2 = 0$ hold.

\item[(ii)] 
If $p \geq 2$ and $q > 2$, then for an almost Grassmann structure 
$AG (p-1, p+q-1)$ to be $\beta$-semiintegrable, it is necessary 
and sufficient that the conditions  $a_\beta = b_\beta^1 = 
b_\beta^2 = 0$ hold. 

\end{description}
\end{theorem}

{\sf Proof.} We prove  part (i) of theorem. The proof 
of part (ii) is similar.

Suppose that $\theta_\alpha, \; \alpha = 1, \ldots, p$, are 
basis forms of the integral 
submanifolds $V_\alpha, \; \dim V_\alpha =p$, 
of the distribution $\Delta_\alpha$ 
appearing in Definition {\bf 1.2}. Then 

\begin{equation}\label{eq:2.9}
\omega_\alpha^i = s^i \theta_\alpha, \;\;\;\;\;
 \alpha = 1, \ldots, p; \;\; i = p+1, \ldots , p+q
\end{equation}
where $\theta_\alpha$ are basis forms on the submanifold $V_\alpha$.

For the structure $AG (p-1, p+q-1)$ to be 
$\alpha$-semiintegrable,  it is necessary and 
sufficient that system (2.9) be completely integrable. 
Taking the exterior derivatives of equations (2.9) by means of 
structure equations (1.4), we find that 

\begin{equation}\label{eq:2.10}
(d s^i + s^j \omega_j^i - s^i \omega) \wedge \theta_\alpha 
+ s^i (d \theta_\alpha - \omega_\alpha^\beta \wedge \theta_\beta) 
= a_{\alpha jk}^{i \beta\gamma} s^j s^k \theta_\beta \wedge 
\theta_\gamma.
\end{equation}
It follows from these equations that 
\begin{equation}\label{eq:2.11}
d \theta_\alpha -  \omega_\alpha^\beta \wedge \theta_\beta 
= \varphi_\alpha^\beta  \wedge \theta_\beta,
\end{equation}
where $\varphi_\alpha^\beta$ is an 1-form that is not 
expressed in terms of  the basis forms $\theta_\alpha$. 

For brevity, we set 
\begin{equation}\label{eq:2.12}
 \varphi^i = d s^i + s^j \omega_j^i - s^i \omega.
\end{equation}
Then the exterior quadratic equation (2.10) takes the form 
\begin{equation}\label{eq:2.13}
(\delta_\alpha^\beta \varphi^i + s^i \varphi_\alpha^\beta)  
\wedge \theta_\beta 
= a_{\alpha jk}^{i \beta\gamma} s^j s^k \theta_\beta \wedge 
\theta_\gamma.
\end{equation}
 From (2.13) it follows that for $\theta_\alpha = 0$, 
the 1-form $\delta_\alpha^\beta \varphi^i + s^i \varphi_\alpha^\beta$ 
vanishes:
\begin{equation}\label{eq:2.14}
\delta_\alpha^\beta \varphi^i (\delta) 
+ s^i \varphi_\alpha^\beta (\delta) = 0.
\end{equation}
Contracting equation (2.14) with respect to the indices 
$\alpha$ and $\beta$, we find that 
\begin{equation}\label{eq:2.15}
\varphi^i = - s^i \varphi (\delta), \;\; 
 \varphi_\alpha^\beta = \delta_\alpha^\beta \varphi (\delta),
\end{equation}
where we set $\varphi (\delta) = 
\frac{1}{p} 
\varphi_\gamma^\gamma (\delta)$.

It follows from (2.15) that on the subvariety $V_\alpha$, 
 the 1-forms  $\varphi^i$ and $\varphi_\alpha^\beta$ 
can be written as follows:
\begin{equation}\label{eq:2.16}
\varphi^i = - s^i \varphi + s^{i\beta} \theta_\beta, 
\;\;  \varphi_\alpha^\beta = \delta_\alpha^\beta \varphi  
+ \widehat{s}_\alpha^{\beta\gamma} \theta_\gamma.
\end{equation}
Substituting these expressions into equations (2.11) and (2.12), 
we find that 
\begin{equation}\label{eq:2.17}
d \theta_\alpha -  \omega_\alpha^\beta \wedge \theta_\beta 
= \varphi  \wedge \theta_\alpha 
+  s_\alpha^{\beta\gamma} \theta_\gamma \wedge \theta_\beta
\end{equation}
where 
$s_\alpha^{\beta\gamma} = \widehat{s}_\alpha^{[\beta\gamma]}$ 
and 
\begin{equation}\label{eq:2.18}
  d s^i + s^j \omega_j^i - s^i \omega 
=  - s^i \varphi  + s^{i\beta} \theta_\beta.
\end{equation}
Substituting (2.17) and (2.18) into 
equation (2.10), we obtain
\begin{equation}\label{eq:2.19}
- s^i s_\alpha^{\beta\gamma} 
- \delta_\alpha^{[\beta} s^{|i|\gamma]} 
= a_{\alpha jk}^{i [\beta\gamma]} s^j s^k.
\end{equation}
Contracting equation (2.19) with respect to 
the indices $\alpha$ and $\beta$, we obtain
$$
- 2 s^i s_\alpha^{\alpha\gamma} - p s^{i\gamma} + 
s^{i\gamma} = 0,
$$
 from which it follows that 
\begin{equation}\label{eq:2.20}
 s^{i\gamma} = s^i s^\gamma,
\end{equation}
where we set 
$s^\gamma = - 
\frac{2}{p-1} s_\alpha^{\alpha\gamma}$.
Substituting (2.20) into (2.19), we find that 
\begin{equation}\label{eq:2.21}
s^i (\delta_\alpha^\gamma s^\beta 
- \delta_\alpha^\beta s^\gamma - 2 s_\alpha^{\beta\gamma}) 
= 2 a_{\alpha jk}^{i [\beta\gamma]} s^j s^k. 
\end{equation}
It follows that 
\begin{equation}\label{eq:2.22}
\delta_\alpha^\gamma s^\beta 
- \delta_\alpha^\beta s^\gamma - 2 s_\alpha^{\beta\gamma} 
=  s_{\alpha j}^{\beta\gamma} s^j, 
\end{equation}
where $s_{\alpha j}^{\beta\gamma} 
= - s_{\alpha j}^{\gamma\beta}$. 
Substituting (2.22) into (2.21), we arrive at the equation 
\begin{equation}\label{eq:2.23}
 s_{\alpha (j}^{\beta\gamma} \delta^i_{k)} = 
a_{\alpha (jk)}^{i\beta\gamma}, 
\end{equation}
where the alternation sign in the right-hand side is dropped 
by Lemma 2.1.

Contracting (2.23) with respect to the indices $i$ and $j$ 
and taking into account of equations (1.9) and (1.10), 
we obtain 
\begin{equation}\label{eq:2.24}
s_{\alpha k}^{\beta\gamma} = 0, 
\end{equation}
from which, by (2.23), it follows that 
\begin{equation}\label{eq:2.25} 
a_{\alpha (jk)}^{i\beta\gamma} = 0. 
\end{equation}

This proves that {\em if an almost Grassmann structure 
$AG (p - 1, p + q - 1)$ is $\alpha$-semiintegrable, then its 
torsion tensor satisfies the condition $(2.25)$,  that is, 
$a_\alpha = 0$.} 

Since, by Theorem 2.4, for $p = 2$ the subtensor 
$a_\alpha = 0$,  condition (2.25) is identically satisfied. 
Hence, while proving sufficiency of this condition for 
$\alpha$-semiintegrability, we must assume that $p > 2$.

Let us return to equations (2.17) and (2.18). Substitute 
into equation (2.18) the values $s^{i\beta}$ taken from 
(2.20) and set 
\begin{equation}\label{eq:2.26}
\widetilde{\varphi} = \varphi - s^\beta \theta_\beta.
\end{equation}
In addition, by (2.24), relations (2.22) imply that 
$$
s_\alpha^{\beta\gamma} = \delta_\alpha^{[\gamma} s^{\beta]}.
$$
Then equations (2.17) and (2.18) take the form 
\begin{equation}\label{eq:2.27}
d \theta_\alpha -  (\omega_\alpha^\beta 
+ \delta_\alpha^\beta \widetilde{\varphi}) \wedge \theta_\beta 
= 0 
\end{equation}
and 
\begin{equation}\label{eq:2.28}
d s^i + s^j \omega_j^i - s^i (\omega - \widetilde{\varphi}) = 0.
\end{equation}
Taking the 
exterior derivatives of (2.28), we obtain the following 
exterior quadratic equation: 
\begin{equation}\label{eq:2.29}
s^i \Phi + b_{jkl}^{i\gamma\delta} s^j s^k s^l \theta_\gamma 
\wedge \theta_\delta = 0,
\end{equation}
where 
$$
\Phi = d \widetilde{\varphi} - \displaystyle \frac{(p+1)q}{p+q} 
s^k \omega_k^\gamma \wedge \theta_\gamma.
$$
Next, taking the exterior derivatives of (2.27), we find that 
\begin{equation}\label{eq:2.30}
 \Phi \wedge \theta_\alpha + b_{\alpha kl}^{\beta\gamma\delta}
 s^k s^l \theta_\beta \wedge \theta_\gamma \wedge \theta_\delta 
= 0.
\end{equation}
Equation (2.29) shows that the 2-form $\Phi$ 
can be written as 
\begin{equation}\label{eq:2.31}
 \Phi = A_{kl}^{\gamma\delta} s^k s^l  \theta_\gamma \wedge \theta_\delta,
\end{equation}
where the coefficients $A_{kl}^{\gamma\delta}$ are symmetric 
with respect to the lower indices and skew-symmetric 
with respect to the upper indices. 
Substituting this value of the form $\Phi$ 
into equations (2.29) and (2.30), we arrive 
at the conditions 
\begin{equation}\label{eq:2.32}
b_{(jkl)}^{i[\gamma\delta]} 
+ \delta^i_{(j} A_{kl)}^{\gamma\delta} = 0
\end{equation}
and 
\begin{equation}\label{eq:2.33}
b_{\alpha (kl)}^{[\beta\gamma\delta]}
 + \delta_\alpha^{[\beta} A_{kl}^{\gamma\delta]} = 0.
\end{equation}
Contracting equation (2.32) with respect to the indices 
$i$ and $j$ and  equation (2.33) with respect to the indices 
$\alpha$ and $\beta$, we find that 
\begin{equation}\label{eq:2.34}
2 (q + 2) A_{kl}^{\gamma\delta} 
+ b_{kli}^{i\gamma\delta} + b_{kil}^{i\gamma\delta} 
+  b_{lik}^{i\gamma\delta} + b_{lki}^{i\gamma\delta}  = 0
\end{equation}
and 
\begin{equation}\label{eq:2.35}
2 (p - 2) A_{kl}^{\gamma\delta} 
+ b_{\alpha kl}^{\gamma\delta\alpha} 
+ b_{\alpha lk}^{\gamma\delta\alpha} 
+ b_{\alpha kl}^{\delta\alpha\gamma} 
+ b_{\alpha lk}^{\delta\alpha\gamma} = 0.
\end{equation}
Note that for $p = 2$ equation (2.33) becomes an 
identity, and we will not obtain equation (2.35). 

If we add equations (2.34) and (2.35) and apply 
condition (1.12), we find that 
\begin{equation}\label{eq:2.36}
A_{kl}^{\gamma\delta} = 0.
\end{equation}
As a result equations (2.32) and (2.33) take the form
\begin{equation}\label{eq:2.37}
b_{(jkl)}^{i[\gamma\delta]} = 0, \;\; 
b_{\alpha (kl)}^{[\beta\gamma\delta]} = 0.
\end{equation}
By Lemma 2.1, conditions (2.37) are equivalent to the 
conditions
\begin{equation}\label{eq:2.38}
b_{(jkl)}^{i\gamma\delta} = 0, \;\; 
b_{\alpha kl}^{[\beta\gamma\delta]} = 0.
\end{equation}
It follows from equations (2.36) and (2.31) that 
\begin{equation}\label{eq:2.39}
d \widetilde{\varphi} = \displaystyle \frac{(p+1)q}{p+q} 
s^k \omega_k^\gamma \wedge \theta_\gamma.
\end{equation}

Finally, taking the exterior derivatives of equations 
(2.39) and applying (2.27), (2.28), and (1.4), 
we obtain the condition
\begin{equation}\label{eq:2.40}
 c_{(ijk)}^{[\alpha\beta\gamma]} = 0.
\end{equation}
These equations will not be trivial only if $p > 2$. 
But, by Lemma 2.2, conditions (2.40) follow from 
integrability conditions (1.15). 

Thus the system of Pfaffian equations (2.9), 
defining integral submanifolds of an $\alpha$-semiintegrable 
almost Grassmann structure, together with Pfaffian equations 
(2.28) and (2.39) following from (2.9), 
is completely integrable if and only if conditions 
(2.25) and (2.38) are satisfied. 
This proves  part (i). 

As we noted in the beginning, the proof of  part (ii) is 
similar. We note only that the 
equations of integral submanifolds $V_\beta, \; \dim V_\beta 
= q$, of the distribution 
$\Delta_\beta$ appearing in Definition 1.6 can be written in 
the form
\begin{equation}\label{eq:2.41}
\omega_\alpha^i = s_\alpha \theta^i, \;\;\;\;\; 
\alpha = 1, \ldots, p; \; i = p+1, \ldots, p+q,
\end{equation}
where the 1-forms $\theta^i$ are linearly independent on 
the submanifold $V_\beta$. \rule{3mm}{3mm}

\vspace*{2mm}

It follows from our previous considerations and  (2.8) that 
\begin{description}
\item[1.] for $p = 2$ we have $b_\alpha^2 = 0$ and 
$c_\alpha = 0$; 
\item[2.] for $q = 2$ we have $b_\beta^1 = 0$ and 
$c_\beta = 0$; 
\item[3.] for $p > 2$ we have  $c_\alpha =0$;  
\item[4.] for $q > 2$ we have $c_\beta = 0$. 
\end{description}

The last two results follow from conditions (1.15) and Lemma 
2.2. 
These results combined with differential 
 equations which the components of $b^1$ and $b^2$ satisfy 
 imply that the tensors $a_\alpha$, $a_\beta$ and the quantities  
$b_\alpha^1, b_\alpha^2,  b_\beta^1, b_\beta^2$ form 
 the following geometric objects:
$$
\begin{array}{lll}
(a_\alpha, b_\alpha^1),  & \!\! (a_\alpha, b_\alpha^2), & \!\! 
S_\alpha = (a_\alpha, b_\alpha^1, b_\alpha^2),    \\ 
&&\\
(a_\beta, b_\beta^1),  &\!\!  (a_\beta, b_\beta^2), & \!\! 
S_\beta =  (a_\beta, b_\beta^1, b_\beta^2),   
\end{array}
$$
which are  subobjects of the second structural object and 
the complete structural object of the almost Grassmann 
structure.

The following theorem gives the conditions of 
$\alpha$- and $\beta$-semiintegrability in the cases when $p = 2$ or  $q = 2$.

\begin{theorem} 
\begin{description}
\item[(i)] If $p = 2$, then the structure subobject $S_\alpha$ 
consists only of the tensor $b_\alpha^1$, and the vanishing of 
this tensor is necessary and 
sufficient for the almost Grassmann structure $AG (1, q + 1)$ 
to be $\alpha$-semiintegrable. 

\item[(ii)] If $q = 2$, then the structure subobject $S_\beta$ 
consists only of the tensor $b_\beta^2$, and the vanishing 
of this tensor is necessary and 
sufficient for the almost Grassmann structure $AG (p-1, p + 1)$ 
$($which is equivalent to the structure $AG (1, p + 1))$ 
to be $\beta$-semiintegrable. 

\item[(iii)] If $p = q = 2$, then the complete structural 
object $S$ consists only of the tensors 
$b_\alpha^1$ and $b_\beta^2$, and the vanishing 
of one of these tensors is necessary and 
sufficient for the almost Grassmann structure $AG (1, 3)$ 
to be $\alpha$- or $\beta$-semiintegrable, respectively. 
\end{description}
\end{theorem}

{\sf Proof.} 
We will prove part (i). As we have already seen, for 
 $p = 2$, the tensor $a_\alpha$ as well as the quantities 
$b_\alpha^2$ and $c_\alpha$ vanish ($a_\alpha = b_\alpha^2 = 
c_\alpha = 0$), and  the object $b_\alpha^1$ becomes a tensor. 
Thus the vanishing of this tensor is necessary and 
sufficient for the almost Grassmann structure $AG (1, q + 1)$ 
to be $\alpha$-semiintegrable. The proof of part (ii) 
is similar. Part (iii) combines the results of (i) 
and (ii). \rule{3mm}{3mm}

We will make two more remarks:

\begin{description}
\item[1.] The tensors 
$b_\alpha^1$ and $b_\beta^2$ are defined by a third-order 
differential neighborhood of the  almost Grassmann structure.

\item[2.] For $p = q = 2$, as was indicated earlier 
(see Subsection {\bf 1.1}), 
the almost Grassmann structure $AG (1, 3)$ is equivalent to the 
conformal $CO (2, 2)$-structure. Thus  we have the 
following decomposition of its 
complete structural object: $S = b_\alpha^1 \dot{+} b_\beta^2$ 
(see \cite{AG96}, \S 5.1). 
\end {description}

Note also that in \cite{D93} and \cite{D94} the author assumed that 
an almost Grassmann structure $AG(p-1, p+q-1)$ is semiintegrable if 
and only if one of two curvature forms $\Omega_i^j$ or 
$\Omega_\alpha^\beta$ 
vanishes. However, our previous considerations as well as 
formula (3.18) below show that this will be the case only if $p = q = 2$.

\section{Existence of  Semiintegrable 
Almost Grassmann 
\newline  
Structures}

\setcounter{equation}{0}

{\bf 1.} We will prove the existence of 
four-dimensional and $(2p)$- and $(2q)$-dimensional semiintegrable 
almost Grassmann structures, where $p, q >2$,  by constructing 
examples of such structures. 

In order to prove that a certain almost Grassmann structure 
is $\alpha$- or $\beta$-semiintegrable, we will use two 
methods. 

The {\em first method} is to check whether conditions 
of $\alpha$- or $\beta$-semiintegrability 
outlined in Theorems 2.5 and 2.6 are satisfied.

When we apply this method,  we will need the following lemma:

\begin{lemma} If the forms $\omega, \omega_\alpha^\beta$, 
and $\omega_j^i$, occurring in equations 
$(1.4)$ are principal forms, that is,
$$
\omega, \omega_\alpha^\beta, \omega_j^i \equiv 0 \pmod{\omega_\alpha^i}, 
$$ 
 then the form  $\omega$ can be 
reduced to $0$, and the forms $\omega_i^\alpha, \alpha = 1, 2, 
i = 3,4$, become principal forms with a symmetric matrix of 
coefficients. 
\end{lemma}

{\sf Proof.} Suppose that $\omega$ is expressed in terms 
of the basis forms $\omega^i_\alpha$:
\begin{equation}\label{eq:3.1}
 \omega = a \omega_1^3 + b \omega_2^3 + c \omega_1^4 
+ e \omega_2^4. 
\end{equation}
If we take exterior derivative of this equation, apply Cartan's 
lemma to the obtained exterior quadratic equation, and 
set $\omega^i_\alpha = 0$, we will obtain the following 
Pfaffian equations:
$$
\renewcommand{\arraystretch}{1.3}
\begin{array}{ll}
\delta a + \pi_3^1 = 0, & \delta b + \pi_3^2 = 0, \\
\delta c + \pi_4^1 = 0, & \delta e + \pi_4^2 = 0, 
\end{array}
\renewcommand{\arraystretch}{1}
$$
where $\delta$ is the symbol of differentiation 
with respect to fiber parameters, and $\pi_i^\alpha 
= \omega_i^\alpha|_{\omega_\gamma^k = 0}$. 

This implies that the quantities $a, b, c$ and $e$ 
can be reduced to 0, $a=b=c=e=0$. Let us show, 
for example, that $a$ can be reduced to 0. In fact, 
since the forms $\pi_3^1, \pi_3^2, \pi_4^1$, and $\pi_4^2$ 
are linearly independent, we can set $\pi_3^2 = \pi_4^1 = \pi_4^2 
= 0$ preserving $\pi_3^1 \neq 0$. Since now  $\pi_3^1$ 
depends on one fiber parameter, for an appropriate choice 
of this parameter we have $\pi_3^1 = \delta t$. By integrating 
the equation $\delta a + \delta t = 0$, we obtain 
$a = - t + C$, where $C$ is a constant. For any value of $C$, 
we can take $t = C$, and as a result, we will get $a = 0$. 

Conversely, if the fiber parameters are chosen in such a way 
that $a = 0$, we get from the above differential 
equations that $\pi_3^1 = 0$. 

The remaining reductions 
$b = c = e = 0$ can be proved in a similar manner. 
This gives $\pi_3^2 = \pi_4^1 = \pi_4^2 = 0$. After these 
reductions equation (3.1) becomes
\begin{equation}\label{eq:3.2}
 \omega = 0.
\end{equation}
Now equations (1.4) and (3.2) imply that 
$$
\omega_i^\alpha \wedge \omega_\alpha^i = 0.
$$
Applying Cartan's 
lemma to this exterior quadratic equation, we obtain 
that the forms $\omega_i^\alpha$ 
become principal forms, and they are expressed in terms 
of the basis forms  $\omega^i_\alpha$ as follows:
\begin{equation}\label{eq:3.3}
\renewcommand{\arraystretch}{1.3}
\left\{
\begin{array}{ll}
\omega_3^1 = A_1 \omega_1^3 + A_2 \omega_2^3 + 
             A_3 \omega_1^4 + A_4 \omega_2^4, \\
\omega_3^2 = A_2 \omega_1^3 + B_2 \omega_2^3 + 
             B_3 \omega_1^4 + B_4 \omega_2^4, \\
\omega_4^1 = A_3 \omega_1^3 + B_3 \omega_2^3 + 
             C_3 \omega_1^4 + C_4 \omega_2^4, \\
\omega_4^2 = A_4 \omega_1^3 + B_4 \omega_2^3 + 
             C_4 \omega_1^4 + E_4 \omega_2^4. \rule{3mm}{3mm}
\end{array}
\right.
\renewcommand{\arraystretch}{1}
\end{equation}

Note that the conditions for fiber forms in Lemma 3.1 mean 
the group $G$ of admissible transformations of first-order frames 
is reduced to the identity group, $G = \{e\}$, and instead 
of a {\em fibration} of first-order frames associated with 
an $AG (1, 3)$-structure, we have a {\em distribution} 
of such frames. In other words, now with any point of 
a manifold $M^4$ on which the $AG (1, 3)$-structure is given 
only {\em one frame} is associated. 

The normalization (3.2) and relations (3.3) 
implied by this normalization means that the subgroup 
 ${\bf T} (4)$ of translations contained in 
the prolonged group $G'$ (see \cite{AG96}, p. 274) 
is also reduced to the identity group, and $G' = G = \{e\}$. 
Furthermore,  only  one  second-order frame is associated with 
any point   $x \in M^4$. 
Moreover, the normalization (3.2)--(3.3) singles out 
a pseudo-Riemannian metric $g$ of signature $(2, 2)$ 
that is concordant with the almost Grassmann structure 
$AG (1, 3)$ given on the manifold $M^4$. 

When we  use the first method, we  also need to use 
relations among the components of the tensor $b = \{b^1, b^2\}$ 
that follow from equations (1.12) and (1.13). There are 16 
equations (1.12) and 256 equations (1.13). In order to 
 get the relations along components of the 
tensor $b$,  first we  
prove the following two lemmas which list 10 independent relations among 
relations (1.12) and 16 independent relations among relations (1.13).

\begin{lemma} For $p = q = 2$, equations $(1.12)$ take the form:
\begin{equation}\label{eq:3.4}
\renewcommand{\arraystretch}{1.5}
\left\{
\begin{array}{ll}
b_{233}^{121} - b_{343}^{411} = 0, & 
b_{133}^{212} - b_{343}^{422} = 0, \\
b_{244}^{121} - b_{434}^{311} = 0, & 
b_{144}^{212} - b_{434}^{322} = 0, 
\end{array}
\right.
\renewcommand{\arraystretch}{1}
\end{equation}
\begin{equation}\label{eq:3.5}
\renewcommand{\arraystretch}{1.5}
\left\{
\begin{array}{ll}
2 b_{133}^{112} - b_{343}^{412} - b_{343}^{421} = 0, & 
2 b_{144}^{112} - b_{434}^{312} - b_{434}^{321} = 0,  \\
2 b_{443}^{411} - b_{234}^{121} - b_{243}^{121} = 0, & 
2 b_{334}^{322} - b_{143}^{212} - b_{134}^{212} = 0,  
\end{array}
\right.
\end{equation}
 and
\begin{equation}\label{eq:3.6}
\renewcommand{\arraystretch}{1.5}
\left\{
\begin{array}{ll}
2 b_{134}^{112} + b_{234}^{122} - 2b_{334}^{312} - b_{344}^{412} 
+ b_{143}^{211}  - b_{433}^{321} = 0,  \\
 b_{134}^{211} + 2b_{234}^{221} - 2b_{334}^{321} - b_{344}^{421} 
 + b_{243}^{122} - b_{433}^{312} = 0.  
\end{array}
\right.
\renewcommand{\arraystretch}{1}
\end{equation} 
\end{lemma}

{\sf Proof.} The proof is straightforward. 
Equations (3.4) are obtained from (1.12) taking 
$\gamma = \delta = 1, k = l = 3; 
\gamma = \delta = 2, k = l = 3; 
\gamma = \delta = 1, k = l = 4$, 
and $\gamma = \delta = 2, k = l = 4$,  respectively. 
Equations (3.5) are obtained from (1.12) taking  
$\gamma = 1, \delta = 2, k = l = 3; 
\gamma = 1, \delta = 2, k = l = 4; \gamma = \delta = 1, 
k = 3, l = 4$, and $\gamma = \delta = 2, k = 3, l = 4$, 
respectively. 
Finally, equations (3.6) are obtained from (1.12) taking 
$\gamma = 1, \delta = 2, k = 3, l = 4$ and 
$\gamma = 2, \delta = 1, k = 3, l = 4$,  respectively. 
The remaining equations (1.12) do not give new conditions. 
Note that while obtaining (3.4)--(3.6), one should use 
conditions (1.10) and (1.11). 

The remaining relations (1.12) are satisfied identically or 
lead to the same relations (3.4)--(3.6).
\rule{3mm}{3mm}

\begin{lemma} For $p = q = 2$, equations $(1.13)$ take the form:
\begin{equation}\label{eq:3.7}
\renewcommand{\arraystretch}{1.5}
\left\{
\begin{array}{ll}
b_{233}^{121} - b_{334}^{411} = 0, & 
b_{133}^{212} - b_{334}^{422} = 0, \\
b_{244}^{121} - b_{443}^{311} = 0, & 
b_{144}^{212} - b_{443}^{322} = 0, 
\end{array}
\right.
\renewcommand{\arraystretch}{1}
\end{equation}
\begin{equation}\label{eq:3.8}
\renewcommand{\arraystretch}{1.5}
\left\{
\begin{array}{ll}
b_{133}^{121} = b_{343}^{412} +  b_{433}^{421}, & 
b_{233}^{212} = b_{343}^{421} +  b_{433}^{412}, \\
b_{144}^{121} = b_{434}^{312} +  b_{344}^{321}, &
b_{244}^{212} = b_{343}^{421} +  b_{434}^{321}, \\
 b_{334}^{311} = b_{234}^{211} +  b_{243}^{112}, & 
b_{443}^{411} = b_{234}^{112} +  b_{243}^{211}, \\
b_{334}^{322} = b_{134}^{122} +  b_{143}^{221}, &
b_{443}^{422} = b_{143}^{122} +  b_{134}^{221}, \\
\end{array}
\right.
\renewcommand{\arraystretch}{1}
\end{equation}
and 
\begin{equation}\label{eq:3.9}
\renewcommand{\arraystretch}{1.5}
\left\{
\begin{array}{ll}
b_{334}^{312} + b_{433}^{321} = b_{234}^{212} + b_{243}^{122}, & 
b_{143}^{211} + b_{134}^{112} = b_{344}^{421} + b_{434}^{412},  \\
b_{334}^{321} + b_{433}^{312} = b_{134}^{121} + b_{143}^{211}, & 
b_{234}^{221} + b_{243}^{122} = b_{344}^{412} + b_{434}^{421}.  
\end{array}
\right.
\end{equation}
 \end{lemma}

{\sf Proof.} The proof is also straightforward. 
Equations (3.7) are obtained from (1.13)  
$\alpha = \gamma = 2, \beta = \delta = 1, i = j = 4, k = l = 3; 
\alpha = \gamma = 1, \beta = \delta = 2, i = j = 4, k = l = 3; 
\alpha = \gamma = 2, \beta = \delta = 1, i = j = 3, k = l = 4;$ 
and 
$\alpha = \gamma = 1, \beta = \delta = 2, i = j = 3, k = l = 4$, 
  respectively. 

Equations (3.8) are obtained from (1.13) taking 
$\alpha = \beta = \delta = 1, \gamma = 2, i = j = 4,  k = l = 3; 
\alpha = \beta = \gamma = 2, \delta = 1, i = k = 4,  j = l = 3; 
\alpha = \beta = \delta = 1, \gamma = 2, i = j = 3,  k = l = 4; 
\alpha = \beta = \gamma = 2, \delta = 1, i = k = 4,  j = l = 3; 
\alpha = \beta = 2, \gamma = \delta = 1, i = j = k = 3, l = 4; 
\alpha = \beta = 2, \gamma = \delta = 1, i = j = k = 4, l = 3; 
\alpha = \beta = 1, \gamma = \delta = 2, i = j = k = 3, l = 4;$ 
and
$\alpha = \beta  = 1, \delta = \gamma = 2, i = j = k = 4, l = 3$, 
respectively. 

Finally, equations (3.9) are obtained from (1.13) taking 
$ \alpha = \beta = \delta = 2, \gamma = 1, i = j = k = 3,  l = 3; 
\alpha = \beta = \gamma = 1, \delta = 2, i = k = l = 4,  k = 3;
\alpha = \beta = \delta = 1, \gamma = 2, i = j = k = 3,  l = 4$; 
and 
$\alpha = \beta = \gamma = 2, \delta = 1, i = j = l = 4,  k = 3$,
respectively. 

The remaining relations (1.13) are satisfied identically or 
lead to the same relations (3.7)--(3.9). 
\rule{3mm}{3mm} 

\begin{theorem} The components of the tensor $b$ of an almost Grassmann structure $AG (1, 3)$ satisfy the following conditions:
\begin{equation}\label{eq:3.10}
b_{233}^{121} = b_{334}^{411} = b_{133}^{212} = b_{334}^{422} 
= b_{244}^{121} = b_{443}^{311} = b_{144}^{212} = b_{443}^{322} = 0, 
\end{equation}
\begin{equation}\label{eq:3.11}
\renewcommand{\arraystretch}{1.5}
\left\{
\begin{array}{ll}
b_{133}^{112} = 0, & b_{343}^{412} = b_{433}^{412} = b_{334}^{412}, \\ 
b_{144}^{112} = 0, & b_{443}^{321} = b_{344}^{321} = b_{434}^{321}, \\ 
b_{443}^{411} = 0, & b_{243}^{211} = b_{243}^{122} = b_{243}^{121}, \\ 
b_{334}^{322} = 0, & b_{134}^{122} = b_{134}^{221} = b_{134}^{212},  
\end{array}
\right.
\renewcommand{\arraystretch}{1}
\end{equation}
and 
\begin{equation}\label{eq:3.12}
\renewcommand{\arraystretch}{1.5}
\left\{
\begin{array}{ll}
b_{134}^{112} = b_{134}^{211} = b_{134}^{121} = b_{243}^{122},\\  
b_{433}^{321} = b_{334}^{321} = b_{343}^{321} = b_{344}^{412}.
\end{array}
\right.
\end{equation}
 \end{theorem}

{\sf Proof.} To prove this theorem, we  combine results 
of Lemmas 3.2 and 3.3. For example, the first two relations 
(3.10) are obtained by comparing the first equations of 
(3.4) and (3.7), and the first two equations (3.11) are obtained 
by comparing the first equation of (3.5) with 
 the sum of first two equations of (3.8).

To prove relations (3.12), we first note that 
it follows from (3.9) and (3.6) that 
\begin{equation}\label{eq:3.13}
\renewcommand{\arraystretch}{1.5}
\left\{
\begin{array}{ll}
2b_{334}^{312} + b_{344}^{412} + b_{433}^{321} = 0, &
2b_{134}^{112} + b_{234}^{122} + b_{143}^{211} = 0, \\
2b_{334}^{321} + b_{344}^{421} + b_{433}^{312} = 0, &
2b_{234}^{221} + b_{143}^{122} + b_{134}^{211} = 0.
\end{array}
\right.
\end{equation}
In fact, by adding the first two equations of (3.9) we find 
that 
$$
2b_{334}^{312} + b_{434}^{412} + b_{433}^{321} = 
- (2b_{134}^{112} + b_{234}^{122} + b_{143}^{211}).
$$
Comparing this equation with the first equation of (3.6), 
we arrived to the first two equations of (3.13). The remaining 
equations of (3.13) are obtained from the last two equations 
of (3.9) and the second equation of (3.6).

If we add the equations of the first and the second column 
of (3.13), we easily find that
\begin{equation}\label{eq:3.14}
b_{334}^{321} = b_{343}^{321}, \;\; b_{134}^{121} = b_{134}^{112}.
\end{equation}
Next adding equations of the first row and the first column of 
(3.9) and taking into account of relations (3.13), we obtain
 \begin{equation}\label{eq:3.15}
 b_{134}^{211} = b_{243}^{122}, \;\;
 b_{433}^{321}= b_{344}^{412}.
\end{equation}
Equations (3.15) and the first two equations (3.13) imply that
\begin{equation}\label{eq:3.16}
b_{134}^{121}  = b_{243}^{122}, \;\;
b_{343}^{321}= b_{344}^{412}.
\end{equation}
Equations (3.15) and (3.16) gives relations 
(3.12). \rule{3mm}{3mm}

It is easy to see from relations (3.10)--(3.12) that 
there are only 10 independent components of the curvature tensor 
$b$ of an almost Grassmann structure $AG (1, 3)$, and that 
we were able to obtain relations (3.10)--(3.12) only because 
for $p = q = 2$ conditions (1.13) become conditions for 
the components of the curvature tensor 
$b$ of an almost Grassmann structure $AG (1, 3)$.

We  make the following remark on a four-dimensional almost Grassmann 
structure $AG (1, 3)$.  Since an $AG (1, 3)$-structure 
is equivalent to a $CO (2, 2)$-structure, we can calculate 
the conformal curvature tensor $C = C_\alpha \dot{+} C_\beta$, where 
$C_\alpha = \{a_0, a_1, a_2, a_3, a_4, a_5\}$ and 
$C_\beta = \{b_0, b_1, b_2, b_3, b_4, b_5\}$ (see \cite{AG96}, Section {\bf 
5.1}) of the latter. 

We impose the following relations between the basis forms 
$\omega^1,  \omega^2, \omega^3, \omega^4$ of a 
$CO (2, 2)$-structure and the basis forms 
$\omega^3_1, \omega^3_2, \omega^4_1,  \omega^4_2$ of an 
equivalent $AG (1, 3)$-structure: 
\begin{equation}\label{eq:3.17}
\renewcommand{\arraystretch}{1.5}
\left\{
\begin{array}{ll}
  \omega^1 = \frac{1}{\sqrt{2}} \omega^3_1, &
\omega^2 = \frac{1}{\sqrt{2}} \omega^3_2,\\
\omega^3 = \frac{1}{\sqrt{2}} \omega^1_4, &
\omega^4 = \frac{1}{\sqrt{2}} \omega^4_2.
\end{array}
 \right.
\renewcommand{\arraystretch}{1}
\end{equation} 
The factor $\frac{1}{\sqrt{2}}$ can be explained by the fact 
that the metric of a $CO (2, 2)$-structure usually is written
in the form $g = 2(\omega^1 \omega^4 - \omega^2 \omega^3$) 
and that the metric of an equivalent $AG (1, 3)$-structure 
has the form 
$g = \omega^3_1 \omega_2^4 - \omega_3^2 \omega_1^4$ 
(see Ch. 5 and 7 of \cite{AG96}), and   relations 
(3.13) make these metrics equal.

 The calculations  
involving the apparatus developed in Ch. 5 of the book 
\cite{AG96} and the formulas (1.4)--(1.5), (3.4)--(3.5), 
give the following  relations between 
the independent 10 components of the curvature 
tensor $b$ of the structure $AG (1, 3)$ (see Theorem {\bf 3.4}) 
and independent 10 components of the tensor $C$ of the 
equivalent pseudoconformal structure $CO (2, 2)$:
\begin{equation}\label{eq:3.18}
\renewcommand{\arraystretch}{1.5}
\left\{
\begin{array}{ll}
a_0 = b_{333}^{412},\\ 
a_1 = b_{433}^{412} = b_{334}^{412} = b_{343}^{412}, \\
a_2 = b_{433}^{321} = b_{334}^{321} = b_{343}^{321} = 
     b_{344}^{412}, \\
a_3 = b_{344}^{321} = b_{434}^{321} = b_{443}^{321}, \\
a_4 =  b_{444}^{321}, \\ 
b_0 = b_{243}^{111}, \\
b_1 = b_{243}^{211} = b_{243}^{112} = b_{243}^{121}, \\
b_2 = b_{134}^{112} = b_{134}^{211} = b_{134}^{121} =  
      b_{243}^{122}, \\
b_3 = b_{134}^{122} = b_{134}^{212} = b_{134}^{221}, \\
b_4 = b_{134}^{222}.
\end{array}
 \right.
\renewcommand{\arraystretch}{1}
\end{equation}

The conditions $C_\alpha = 0, C_\beta = 0$ and $C = 0$ 
are necessary and sufficient for a $CO (2, 2)$-structure 
(and consequently for an equivalent ($AG (1, 3)$-structure) to be 
$\alpha$-semiintegrable, $\beta$-semiintegrable and locally flat, 
respectively. Moreover, if we find components of 
$C_\alpha$ and $C_\beta$, then by investigating the roots 
of polynomials
$$
C_\alpha (\lambda) = a_0 \lambda^4 - 4a_1 \lambda^3 + 6 a_2 \lambda_2 - 4a_3 
\lambda + a_4
$$
and 
$$
C_\beta (\mu) = b_0 \lambda^4 - 4b_1 \lambda^3 + 6 b_2 \lambda_2 
- 4b_3 \lambda + b_4
$$
 we can make some additional conclusions 
on integrability of the  distributions $\Delta_\alpha (\lambda)$ 
and $\Delta_\beta (\mu)$ of the fiber bundles $E_\alpha$ and $E_\beta$ 
associated with a $CO (2, 2)$-structure (see \cite{AG96}, Section {\bf 5.4}).

Note that equations (3.18) show that for a four-dimensional almost 
Grassmann structure $AG (1, 3)$ and for an equivalent pseudoconformal 
structure $CO (2, 2)$ we have

\begin{description}
\item[i)] The condition of $\alpha$-semiintegrability 
$C_\alpha = 0$ is equivalent to the condition 
$b_{jkl}^{i\beta\gamma} = 0$, that is, $\Omega_i^j = 0$.

\item[ii)] The condition of $\beta$-semiintegrability 
$C_\beta = 0$ is equivalent to the condition 
$b_{\alpha jk}^{\beta\gamma\delta} = 0$, that is, 
 $\Omega_\alpha^\beta = 0$.

\item[iii)] The condition of local flatness $C = 0$ is equivalent 
to the vanishing of the tensor $b = \{b^1_\alpha, b_\beta^2\}$ 
(cf. the end of Section {\bf 1}). 
\end{description}

Note that only for a four-dimensional almost 
Grassmann structure $AG (1, 3)$ the conditions of 
$\alpha$- and $\beta$-semiintegrability are equivalent to 
$\Omega_i^j = 0$ and  $\Omega_\alpha^\beta = 0$, respectively: 
 only for such structures the components of the tensor 
$b$ satisfy relations (3.10)--(3.12) by means of which 
the conditions $b_{\alpha kl}^{[\beta\gamma\delta])} = 0$ 
and $b^{i \gamma\delta}_{[jkl]} = 0$ of $\alpha$- and 
$\beta$-semiintegrability imply the vanishing the components 
$b_{\alpha kl}^{\beta\gamma\delta}$ 
and $b^{i \gamma\delta}_{jkl}$ themselves.

{\bf 2.} The {\em second method} is the method of direct 
integration of equations (2.9) and (2.41) 
defining integral submanifolds of distributions $\Delta_\alpha$ 
and $\Delta_\beta$ of semiintegrable 
almost Grassmann structures. Let us describe this method 
in more detail.

In the proof of Theorem 2.5 we wrote the equations 
of the submanifolds $V_\alpha$ and $V_\beta$ in the form 
(2.9) and (2.41), respectively. Note that $\dim V_\alpha = p$ and 
$\dim V_\beta = q$.  
An almost Grassmann structure is $\alpha$-semiintegrable (or 
$\beta$-semiintegrable) if and only if the system of equations 
(2.9) (respectively, (2.41)) is completely integrable.

We write the matrix of basis forms of the almost Grassmann 
manifold $AG (p - 1, p + q -1)$ in more detail:
\begin{equation}\label{eq:3.19}
\renewcommand{\arraystretch}{1.5}
(\omega_\alpha^i) = \pmatrix{
\omega_1^{p+1} & \omega_2^{p+1} & \ldots & \omega_p^{p+1} \cr 
\omega_1^{p+2} & \omega_2^{p+2} & \ldots & \omega_p^{p+2} \cr 
\multispan4\dotfill \cr
\omega_1^{p+q} & \omega_2^{p+q} & \ldots & \omega_p^{p+q}  
}.
\renewcommand{\arraystretch}{1}
\end{equation} 
Here $\alpha = 1, \ldots , p$ is the column number, and 
$i = p+1, \ldots , p+q$ is the row number.

The condition (2.9) of $\alpha$-integrability means that 
on integral submanifolds $V_\alpha$ of the distribution 
$\Delta_\alpha$ (see Definition 1.2) the rows of 
the matrix $(\omega_\alpha^i)$ are proportional, and 
the entries of every nonzero row are basis forms on 
$V_\alpha$. 

The condition (2.41) of $\beta$-integrability means that 
on integral submanifolds $V_\beta$ of the distribution 
$\Delta_\beta$ (see Definition 1.2) the columns of 
the matrix $(\omega_\alpha^i)$ are proportional, and 
the entries of every nonzero row are basis forms on 
$V_\beta$. 

For $p = q = 2$, equations (2.9) of submanifolds $V_\alpha$ 
can be written in the form 
\begin{equation}\label{eq:3.20}
\lambda \omega_1^3 + \omega_1^4 = 0, \;\; 
\lambda \omega_2^3 + \omega_2^4 = 0,
\end{equation} 
where $\lambda = -\frac{s^4}{s^3}$ and $\omega_1^3 \wedge 
\omega_2^3 \neq 0$, and equations (2.41) of submanifolds $V_\beta$ 
can be written in the form 
\begin{equation}\label{eq:3.21}
\mu \omega_2^3 + \omega_1^3 = 0, \;\; 
\mu \omega_2^4 + \omega_1^4 = 0,
\end{equation} 
where $\mu = -\frac{s_1}{s_2}$ and $\omega_2^3 \wedge 
\omega_2^4 \neq 0$. 

If $p > 2$ or $q > 2$, the systems (3.20) and (3.21) have 
different forms. 
For example, let us consider the case $p = 2$ and $q = 3$. 
In this case equations (2.9) can be  
written as follows:
\begin{equation}\label{eq:3.22}
\renewcommand{\arraystretch}{1.5}
\left\{
\begin{array}{ll}
\lambda_1 \omega_1^3 + \omega_1^4 = 0, &
\lambda_2 \omega_1^3 + \omega_1^5 = 0, \\
\lambda_1 \omega_2^3 + \omega_2^4 = 0, &
\lambda_2 \omega_2^3 + \omega_2^5 = 0,
\end{array}
\right.
\renewcommand{\arraystretch}{1}
\end{equation} 
where $\lambda_1 = -\frac{s^4}{s^3}$ and 
$\lambda_2 = -\frac{s^5}{s^3}$ and $\omega_1^3 \wedge 
\omega_2^3 \neq 0$, and equations (2.41) 
take the form: 
\begin{equation}\label{eq:3.23}
\renewcommand{\arraystretch}{1.5}
\left\{
\begin{array}{ll}
\mu \omega_2^3 + \omega_1^3 = 0, \\
\mu \omega_2^4 + \omega_1^4 = 0, \\
\mu \omega_2^5 + \omega_1^5 = 0, 
\end{array}
\right.
\renewcommand{\arraystretch}{1}
\end{equation} 
where $\mu = -\frac{s_1}{s_2}$ and $\omega_2^3 \wedge 
\omega_2^4  \wedge \omega_2^5 \neq 0$.

For these cases to prove that an almost Grassmann structure 
is $\alpha$-semiintegrable (resp. $\beta$-semiintegrable) we must 
prove that the system (3.20) or (3.22) (resp. the system (3.21) 
or (3.23)) is completely integrable. If it is possible, we   
integrate these systems and find $\lambda$ or $\lambda_1$ and 
$\lambda_2$ (resp. $\mu$) 
and  closed form equations of submanifolds $V_\alpha$ 
(resp. $V_\beta$). 

{\bf 3.}  We next construct  examples of semiintegrable and 
integrable almost Grassmann structures $AG (1, 3)$. To prove 
that they are semiintegrable, we will apply one of two methods 
indicated above.

\examp{\label{examp:3.5} Suppose that $x, y, u$, and $v$ 
are coordinates in $M^4$, and that 
the basis 1-forms $\omega_\alpha^i$ of an almost Grassmann 
structure $AG (1, 3)$ are
\begin{equation}\label{eq:3.24}
\renewcommand{\arraystretch}{1.3}
\left\{
\begin{array}{ll}
\omega_1^3 = dx + f(u) dy, & \omega_2^3 = dy,\\ 
\omega_1^4 = du, & \omega_2^4 = dv. 
\end{array}
\right.
\renewcommand{\arraystretch}{1}
\end{equation}
Taking exterior derivatives of equations (3.22) by means 
of (1.4), (1.5) and (3.22), we arrive at the following exterior 
quadratic equations:
\begin{equation}\label{eq:3.25}
\renewcommand{\arraystretch}{1.3}
\left\{
\begin{array}{ll}
(\omega +  \omega_1^1 - \omega_3^3)\wedge \omega_1^3  
+ \omega_1^2 \wedge  \omega_2^3 + \omega_1^4 \wedge  \omega_4^3 
= f'(u) \omega_1^4 \wedge  \omega_2^4, \\ 
(\omega - \omega_1^1 - \omega_3^3)\wedge \omega_2^3  
+ \omega_2^1 \wedge  \omega_1^3 + \omega_2^4 \wedge  \omega_4^3 = 0, \\ 
(\omega + \omega_1^1 + \omega_3^3)\wedge \omega_1^4  
+ \omega_1^2 \wedge  \omega_2^4 + \omega_1^3 \wedge  \omega_3^4 = 0, \\ 
(\omega - \omega_1^1 + \omega_3^3)\wedge \omega_2^4  
+ \omega_2^1 \wedge  \omega_1^4 + \omega_2^3 \wedge  \omega_3^4 = 0.  
\end{array}
\right.
\renewcommand{\arraystretch}{1}
\end{equation}
First, equations (3.25) prove that the form $\omega$ is a principal form. Thus, by Lemma 3.1 we have equations (3.2) and (3.3). 
Second, it follows from (3.25), that the forms 
$\omega_1^2, \omega_2^1, \omega_3^4, \omega_4^3, \newline 
\omega_1^1 + \omega_3^3$ and $\omega_1^1 - \omega_3^3$ are 
principal forms. We will write their expressions as follows:
\begin{equation}\label{eq:3.26}
\renewcommand{\arraystretch}{1.3}
\left\{
\begin{array}{cl}
\omega_1^2& = \alpha_1 \omega_1^3 + \alpha_2 \omega_2^3 + 
             \alpha_3 \omega_1^4 + \alpha_4 \omega_2^4, \\
\omega_2^1& = \beta_1 \omega_1^3 + \beta_2 \omega_2^3 + 
             \beta_3 \omega_1^4 + \beta_4 \omega_2^4, \\
\omega_3^4 &= \gamma_1 \omega_1^3 + \gamma_2 \omega_2^3 + 
             \gamma_3 \omega_1^4 + \gamma_4 \omega_2^4, \\
\omega_4^3 &= \delta_1 \omega_1^3 + \delta_2 \omega_2^3 + 
             \delta_3 \omega_1^4 + \delta_4 \omega_2^4, \\ 
\omega_1^1 + \omega_3^3 &
           = \sigma_1 \omega_1^3 + \sigma_2 \omega_2^3 + 
             \sigma_3 \omega_1^4 + \sigma_4 \omega_2^4, \\ 
\omega_1^1 - \omega_3^3 &
           = \tau_1 \omega_1^3 + \tau_2 \omega_2^3 + 
             \tau_3 \omega_1^4 + \tau_4 \omega_2^4. 
\end{array}
\right.
\renewcommand{\arraystretch}{1}
\end{equation}
Substituting (3.2) and (3.26) into equations (3.25) and equating 
coefficients in the independent exterior forms $\omega_\alpha^i 
\wedge \omega_\beta^j$ to 0, we obtain: 
$\alpha_1 =  \alpha_2 =  \alpha_4 = 0;  
\beta_1 =  \beta_2 =  \beta_3 = \beta_4 = 0;  
\gamma_1 =  \gamma_2 =  \gamma_3 = \gamma_4 = 0;  
\delta_1 =  \delta_3 =  \delta_4 = 0;  
\sigma_1 =  \sigma_2 =  \sigma_4 = 0;
 \tau_1 =  \tau_2 =  \tau_3 = \tau_4 = 0$, and $\alpha_3 = 
\delta_2 = \sigma_4 = \frac{1}{2} f' (u)$.  
 As a result, equations (3.26) become
\begin{equation}\label{eq:3.27}
\renewcommand{\arraystretch}{1.3}
\left\{
\begin{array}{ll}
\omega_1^2 = \frac{1}{2} f' (u) \omega_1^4, 
\;\;\;\;\omega_2^1 = 0, \\
\omega_4^3 = \frac{1}{2} f' (u) \omega_2^3, 
\;\;\;\;\omega_3^4 = 0, \\
\omega_1^1 = - \omega_2^2 
= \omega_3^3 = - \omega_4^4 
           = \frac{1}{4} f' (u) \omega_2^4. 
\end{array}
\right.
\renewcommand{\arraystretch}{1}
\end{equation}

Taking exterior derivatives of (3.27) by means of 
(1.4), (3.21) and (3.27), we arrive at the following system of 
exterior quadratic equations:
\begin{equation}\label{eq:3.28}
\renewcommand{\arraystretch}{1.3}
\left\{
\begin{array}{ll}
-\omega_1^3 \wedge  \omega_3^2 - \omega_1^4 \wedge \omega_4^2  
+ \Omega_1^2 = \frac{1}{4} (f'(u))^2 \omega_1^4 
\wedge  \omega_2^4, \\ 
-\omega_2^3 \wedge  \omega_3^1 - \omega_1^4 \wedge \omega_4^2  
+ \Omega_2^1 = 0, \\ 
-\omega_3^1 \wedge  \omega_1^4 - \omega_3^2 \wedge \omega_2^4  
+ \Omega_4^3 = 0, \\ 
-\omega_4^1 \wedge  \omega_1^3 - \omega_4^2 \wedge \omega_2^3  
+ \Omega_4^3 = \frac{1}{2} f''(u) \omega_1^4 \wedge  \omega_2^3 
+ \frac{1}{4} (f'(u))^2 \omega_2^4 \wedge  \omega_2^3, \\ 
-\omega_1^3 \wedge  \omega_3^1 - \omega_1^4 \wedge \omega_4^1  
+ \Omega_1^1 = \frac{1}{4} f''(u) \omega_1^4 \wedge  \omega_2^4, \\ 
-\omega_3^1 \wedge  \omega_1^3 - \omega_4^1 \wedge \omega_1^4  
+ \Omega_3^3 = \frac{1}{4} f''(u) \omega_1^4 \wedge  \omega_2^4,  
\end{array}
\right.
\renewcommand{\arraystretch}{1}
\end{equation}
where $\Omega_\alpha^\beta$ and $\Omega_i^j$ are the curvature 
2-forms defined by (1.5). Substituting the values of these 
2-forms from (1.5) and the values of the forms $\omega_i^\alpha$ 
from (3.3) into equations (3.28) and equating 
coefficients in independent exterior forms $\omega_\alpha^i 
\wedge \omega_\beta^j$ to 0, in addition to 
equations (3.10)--(3.12) (which hold for any 
$AG (1, 3)$-structure) we obtain 
 the following additional conditions:
$$
\renewcommand{\arraystretch}{1.5}
\left\{
\begin{array}{ll}
  b_{243}^{111} = b_{243}^{211} = b_{243}^{112} = b_{243}^{121} 
= b_{134}^{112} = b_{134}^{211} = b_{134}^{121} = b_{243}^{122} =0,\\
b_{134}^{122} = b_{134}^{212} = b_{134}^{221} = b_{134}^{222} 
= b_{333}^{412} = b_{433}^{412} = b_{334}^{412} = b_{343}^{412} =0,\\
b_{433}^{321} = b_{334}^{321} = b_{343}^{321} = b_{344}^{412} 
= b_{444}^{321} = 0, \\
  b_{344}^{321} = b_{434}^{321} = b_{443}^{321} 
     = - \frac{1}{4} f''(u).
\end{array}
\right.
\renewcommand{\arraystretch}{1}
$$
that is, the only nonvanishing components of the object $b 
= \{b^1, b^2\}$ are the components $b_{434}^{312}, 
b_{434}^{321}$ and $b_{344}^{312} = - b_{444}^{412}$. 

In addition, we find that  
all coefficients of (3.3), except $E_4$ and $C_4$, equal 
to 0, and for the coefficients  $E_4$ and $C_4$
 we find the following values:
$$
C_4 = - \frac{1}{4} f''(u), \;\; E_4 = - \frac{1}{4} (f'(u))^2.
$$
As a result, equations (3.3) become 
\begin{equation}\label{eq:3.29}
\renewcommand{\arraystretch}{1.3}
\left\{
\begin{array}{ll}
\omega_3^1 = 0, & \omega_4^1 =  - \frac{1}{4} f''(u) \omega_2^4\\
\omega_3^2 = 0, & \omega_4^2 = - \frac{1}{4} f''(u) \omega_1^4 
- \frac{1}{4} (f'(u))^2 \omega_2^4,
\end{array}
\right.
\renewcommand{\arraystretch}{1}
\end{equation}
and the curvature 2-forms (1.7) become
\begin{equation}\label{eq:3.30}
\renewcommand{\arraystretch}{1.3}
\left\{
\begin{array}{ll}
\Omega_1^1 = \Omega_2^2 =  \Omega_1^2 = \Omega_2^1  
= \Omega_3^4 = 0, \\
 \Omega_3^3 = - \Omega_4^4 = \frac{1}{4} f''(u) 
\omega_1^4 \wedge \omega_2^4, \\
 \Omega_4^3 =  \frac{1}{4} f''(u) \omega_1^3 \wedge \omega_2^4 
-  \frac{1}{4} f''(u) \omega_2^3 \wedge \omega_1^4. 
\end{array}
\right.
\renewcommand{\arraystretch}{1}
\end{equation}
Since $b_{\alpha kl}^{\beta \gamma \delta} = 0$, by (2.41), 
we have  $b_\beta^2 = 0$, and the almost Grassmann structure 
$AG (1, 3)$ under consideration is $\beta$-semiintegrable. 
This structure is not locally flat: in fact, $b_{(333)}^{312} = 
b_{(434)}^{312} = \frac{1}{4} f'' (u) \neq 0$ for a general 
function $f (u)$, and thus $b_\alpha^1 \neq 0$, that is, the 
structure is not $\alpha$-semiintegrable. Note that $f'' (u) = 0$ 
if and only if $f (u) = au + b$, where $a$ and $b$ are constants. 
In this case the structure is $\beta$-semiintegrable, and 
consequently it is locally flat.

As we noted earlier,  a structure 
$AG (1, 3)$ is equivalent to a conformal $CO (2, 2)$-structure.
Note that the $CO (2, 2)$-structure corresponding to the 
$\beta$-semiintegrable structure 
$AG (1, 3)$ we have constructed in this example  is self-dual. 

By proving the local existence of a $\beta$-semiintegrable 
structure $AG (1, 3)$, we also proved a local existence 
of a self-dual $CO (2, 2)$-structure. On the global existence 
of four-dimensional semiintegrable smooth compact oriented 
Riemannian manifolds see \cite{T92} and \cite{LB95}.

We will find now the metric of this $CO (2, 2)$-structure. 

To this end, we recall that the equation of the Segre cone 
of the $AG (1,3)$-structure (or the asymptotic cone of 
the corresponding $CO (2, 2)$-structure) has the form 
(1.2) or $\omega_1^3 \omega_2^4 - \omega_1^4 \omega_2^3 = 0$. 
Thus the fundamental form of  the $CO (2, 2)$-structure 
is 
$$
g = \omega_1^3 \omega_2^4 - \omega_1^4 \omega_2^3,
$$
or by (3.14),
\begin{equation}\label{eq:3.31}
g = dx dv + f (u) dy dv - dy du.
\end{equation}
The quadratic form (3.31) defines on the manifold $M^4$ 
a pseudo-Riemannian metric of signature $(2, 2)$ that 
is conformally equivalent to the almost Grassmann structure 
$AG (1, 3)$ with the basis forms (3.24).

Let us apply the second method to prove that the 
structure (3.24) is $\beta$-semiintegrable and to find 
conditions under which this structure is $\alpha$-semiintegrable. 

For the structure (3.24),  the equations (3.21) take the form 
\begin{equation}\label{eq:3.32}
\mu dy + dx + f(u) dy = 0,  \;\; \mu dv + du = 0,
\end{equation}
where $dy \wedge dv \neq 0$. Exterior differentiation 
of (3.32) gives
\begin{equation}\label{eq:3.33}
(d\mu + f' (u)du) \wedge dy = 0,  \;\; 
(d\mu + f' (u)du) \wedge dv = 0.
\end{equation}
Since $dy \wedge dv \neq 0$, it follows from (3.33) that 
\begin{equation}\label{eq:3.34}
d\mu + f' (u)du = 0.
\end{equation}
The solution of (3.34) is 
\begin{equation}\label{eq:3.35}
\mu = C_1 - f (u),
\end{equation}
where $C_1$ is an arbitrary constant.
Substituting this value of $\mu$ into equations (3.32), 
we find the following differential equations of 
submanifolds $V_\beta$:
\begin{equation}\label{eq:3.36}
C_1 dy + dx = 0,  \;\; (C_1 - f (u)) dv + du = 0.
\end{equation}
By integration we find the following {\em closed form 
equations of submanifolds $V_\beta$}:
\begin{equation}\label{eq:3.37}
x = C_2-C_1 y,  \;\; v = C_3 - \displaystyle 
\int \frac{du}{C_1 - f (u)},
\end{equation}
where $C_2$ and $C_3$ are arbitrary constants.

Thus we proved again the structure $AG (1, 3)$ with 
the basis forms (3.14) is $\beta$-semiintegrable. 
In addition, we proved that {\em the integral submanifolds 
$V_\beta$ of this structure are defined by equations $(3.37)$, 
  the family of these submanifolds  depends 
on three parameters $C_1, C_2$ and $C_3$, and 
through any point $(x_0, y_0, u_0, v_0)$ of $M^4$ there passes an 
one-parameter family of submanifolds $V_\beta$.}

To find conditions of $\alpha$-semiintegrability of this 
structure, we first specialize equations (3.20) for it:
\begin{equation}\label{eq:3.38}
\lambda (dx + f(u) dy) + du = 0,  \;\; \lambda dy + dv = 0,
\end{equation}
where $dx \wedge dy \neq 0$. Exterior differentiation 
of (3.38) gives
\begin{equation}\label{eq:3.39}
(d\lambda + \lambda^2 f' (u) dy) \wedge dx = 0,  \;\; 
(d\lambda + \lambda^2 f' (u) dy) \wedge dy = 0.
\end{equation}
Since $dx \wedge dy \neq 0$, it follows from (3.39) that 
\begin{equation}\label{eq:3.40}
d\lambda + \lambda^2 f' (u) dy = 0.
\end{equation}
Exterior differentiation of (3.40) leads to the following 
equation:
\begin{equation}\label{eq:3.41}
\lambda^3 f'' (u) dx \wedge dy = 0.
\end{equation}
Since $dx \wedge dy \neq 0$, in general the structure (3.24) is 
not $\alpha$-semiintegrable. It will be $\alpha$-semiintegrable 
if and only if $f'' (u) = 0$, or 
\begin{equation}\label{eq:3.42}
f (u) = au + b. 
\end{equation}
Substituting $f (u) = au + b$ into equation (3.40), 
we find that 
\begin{equation}\label{eq:3.43}
d\lambda + \lambda^2 a dy = 0.
\end{equation}
The solution of (3.43) is 
\begin{equation}\label{eq:3.44}
\lambda =   \frac{1}{ay+C_4},
\end{equation}
where $C_4$ is an arbitrary constant. Substituting 
$f (u) = au + b$  and $\lambda$ from (3.44) into 
equations (3.38), we find the following differential equations of 
submanifolds $V_\alpha$:
\begin{equation}\label{eq:3.45}
dx + a d(uy) + b dy + C_4 du = 0,  \;\;  dy + (ay + C_4) dv = 0.
\end{equation}
By integration we find the following {\em closed form 
equations of submanifolds $V_\alpha$}:
\begin{equation}\label{eq:3.46}
x + auy + by + C_4 y = C_5,  \;\; \frac{1}{a} \log (ay + C_4) 
+ v = C_6,
\end{equation}
where $C_5$ and $C_6$ are arbitrary constants.

Thus, we proved again that  in general the structure $AG (1, 3)$ 
with the basis forms (3.24) is 
not $\alpha$-semiintegrable, and that it will be 
$\alpha$-semiintegrable and consequently integrable 
if and only if $f (u) = au + b$. If it is the case, 
{\em the closed form equations 
of  submanifolds  $V_\alpha$  and $V_\beta$ are 
equations $(3.46)$ and $(3.37)$. 
If $f (u) = au + b$,  the  family of submanifolds $V_\alpha$ on 
the manifold $M^4$ carrying the $AG (1, 3)$-structure 
with the basis forms $(3.24)$ depends 
on three parameters $C_4, C_5$ and $C_6$, and 
through any point $(x_0, y_0, u_0, v_0)$ 
of $M^4$ there passes a one-parameter family 
of submanifolds $V_\alpha$.}
}

Let us indicate   three more examples of almost Grassmann 
structures similar to the structure (3.24). 

\examp{\label{examp:3.6} 
\begin{equation}\label{eq:3.47}
\renewcommand{\arraystretch}{1.3}
\left\{
\begin{array}{ll}
\omega_1^3 = dx, & \omega_2^3 = dy,\\ 
\omega_1^4 = du + g (x) dv, & \omega_2^4 = dv. 
\end{array}
\right.
\renewcommand{\arraystretch}{1}
\end{equation}
}
\examp{\label{examp:3.7}
\begin{equation}\label{eq:3.48}
\renewcommand{\arraystretch}{1.3}
\left\{
\begin{array}{ll}
\omega_1^3 = dx, &\omega_2^3 = h(v) dx + dy,\\ 
 \omega_1^4 = du, & \omega_2^4 = dv. 
\end{array}
\right.
\renewcommand{\arraystretch}{1}
\end{equation}
}
and 
\examp{\label{examp:3.8} 
\begin{equation}\label{eq:3.49}
\renewcommand{\arraystretch}{1.3}
\left\{
\begin{array}{ll}
\omega_1^3 = dx, & \omega_2^3 = dy,\\ 
\omega_1^4 = du, & \omega_2^4 = k (y) du + dv. 
\end{array}
\right.
\renewcommand{\arraystretch}{1}
\end{equation}
}
As was the case for  the structure (3.24), {\em each of 
the structures $(3.47)$--$(3.49)$ is $\beta$-semiintegrable and 
in general  not $\alpha$-semiintegrable. These structures will be 
$\alpha$-semiintegrable and consequently integrable 
if and only if the functions $g (x), h(v)$ and $k (y)$ are linear 
functions of $x, v$, and $y$, respectively.}

If we apply formulas (3.18) to Examples {\bf 3.5--3.8}, 
then, we obtain that 
$C_\beta = 0$ for all these examples (that is, the corresponding 
structures are $\beta$-semiintegrable), and that the only 
nonvanishing components of the tensor $C_\alpha$ are 
$a_3$ for Examples {\bf 3.5} ($a_3 = -\frac{f'' (u)}{4}$) 
and {\bf 3.7} ($a_3 = \frac{h'' (v)}{4}$) 
and $a_1$ for Examples {\bf 3.6} ($a_1 = \frac{g'' (x)}{4}$) 
and {\bf 3.8} ($a_1 = -\frac{k'' (y)}{4}$). It follows that 
in the first case the polynomial 
$C_\alpha$ has the triple root $\lambda_1 = \lambda_2 = \lambda_3 
= \infty$ and the  simple root $\lambda_4 = 0$, and in the second 
case it has the simple root $\lambda_1 = \infty$  and the triple 
root $\lambda_2 = \lambda_3 = \lambda_4 = 0$. 
According to Section {\bf 5.4} of the book \cite{AG96}, 
this means that the fiber bundle $E_\alpha$ has two invariant  
distributions $\Delta_\alpha (\infty)$ and $\Delta_\alpha (0)$, 
and the distribution corresponding to a multiple root is 
integrable. Moreover, it is easy to prove that the distribution 
corresponding to a simple root is also integrable.  
For all cases the integral submanifolds $V_\alpha$ are defined by 
the equations $u = C_3, v = C_4$ (for $\lambda = \infty$) 
and by the equations $x = C_1, y = C_2$ (for $\lambda = 0$).

The next example was considered in 
\cite{D93}.  However, since results obtained in 
\cite{D93} contain inaccuracies and the result is wrong 
(according to \cite{D93}, the structure of this example 
is of general type), we  give here a complete 
investigation of this example.

\examp{\label{examp:3.9} Suppose again that $x, y, u$ and $v$ 
are coordinates in $M^4$, and that 
the basis 1-forms $\omega_\alpha^i$ of an almost Grassmann 
structure $AG (1, 3)$ are
\begin{equation}\label{eq:3.50}
\renewcommand{\arraystretch}{1.3}
\left\{
\begin{array}{ll}
\omega_1^3 = dx + f (u) dy, & \omega_2^3 = dy,\\ 
\omega_1^4 = du + g (y) dv, & \omega_2^4 = dv. 
\end{array}
\right.
\renewcommand{\arraystretch}{1}
\end{equation}

The metric of the $CO (2, 2)$-structure which is equivalent 
to the $AG (1, 3)$-structure with basis forms (3.40) is 
\begin{equation}\label{eq:3.51}
g = dx dv - dy du + (f (u) - g (y)) dy dv.
\end{equation}

In this case we have 
\begin{equation}\label{eq:3.52}
\renewcommand{\arraystretch}{1.3}
\left\{
\begin{array}{ll}
\omega_1^2 =   g' (y) \omega_2^3 + \frac{1}{2} f'(u) \omega_1^4 
- \frac{1}{2} f'(u) g (y) \omega_2^4, &
 \\
\omega_4^3 = \frac{1}{2} f'(u) \omega_2^3, \;\;\;\; 
\omega_2^1 = 0, 
\;\;\;\;\omega_3^4 =  0,\\
\omega_1^1 = - \omega_2^2 
= \omega_3^3 = - \omega_4^4  = \frac{1}{4} f' (u) \omega_2^4; 
\end{array}
\right.
\renewcommand{\arraystretch}{1}
\end{equation}
\begin{equation}\label{eq:3.53}
C_4 = - \frac{1}{4} f''(u), \;\; E_4 = \frac{1}{2} f''(u) g (y) 
- \frac{1}{4} (f'(u))^2;
\end{equation}
\begin{equation}\label{eq:3.54}
b_{434}^{312} = - b_{434}^{321}=   \frac{1}{4} f''(u), \;\;
b_{344}^{312} = - b_{444}^{412} =  \frac{1}{4} f''(u);
\end{equation}
\begin{equation}\label{eq:3.55}
\renewcommand{\arraystretch}{1.3}
\left\{
\begin{array}{ll}
\omega_3^1 = 0,  &\omega_4^1 = 
 -\frac{1}{4} f''(u) \omega_2^4,  \\
\omega_3^2 = 0,   & 
\omega_4^2 = - \frac{1}{4} f''(u) \omega_1^4 
+ \Bigl[\frac{1}{2} f''(u) g (y) - \frac{1}{4} (f'(u))^2\Bigr] 
\omega_2^4;
\end{array}
\right.
\renewcommand{\arraystretch}{1}
\end{equation}
\begin{equation}\label{eq:3.56}
\renewcommand{\arraystretch}{1.3}
\left\{
\begin{array}{ll}
\Omega_1^1 = \Omega_2^2 =  \Omega_1^2 = \Omega_2^1  
= \Omega_3^4 = 0, \\
 \Omega_3^3 = - \Omega_4^4 = \frac{1}{4} f''(u) 
\omega_1^4 \wedge \omega_2^4, \\
 \Omega_4^3 = \frac{1}{4} f''(u) \omega_1^3 \wedge \omega_2^4 
-  \frac{1}{4} f''(u) \omega_2^3 \wedge \omega_1^4. 
\end{array}
\right.
\renewcommand{\arraystretch}{1}
\end{equation}

We  have $b_{\alpha kl}^{\beta \gamma \delta} = 0$, and by 
(2.41), $b_\alpha^1 = 0$, and $ - b_{444}^{412} = b_{(344)}^{312} 
= -\frac{1}{4} f''(u) \neq 0$. This implies that $b_\beta^2 
\neq 0$, and {\em the almost Grassmann structure $AG (1, 3)$ with 
the basis forms $(3.50)$ is $\beta$-semiintegrable but not 
locally flat.} 

If we apply the  method of direct integration, we find 
that equations (3.21) take the form 
$$
\mu dy + dx + f(u) dy = 0,  \;\; \mu dv + du  + g(y) dv = 0,
$$
where  $dy \wedge dv \neq 0.$ Exterior differentiation 
of these equations leads to 
$$
(d\mu + f' (u)du) \wedge dy = 0,  \;\; 
(d\mu + g' (y)dy) \wedge dv = 0.
$$
It follows that 
$$
d\mu + f' (u) du + g' (y) dy = 0.
$$
By integration of the first equation, we find that 
\begin{equation}\label{eq:3.57}
\mu = C_1 - f (u) - g(y),
\end{equation}
where $C_1$ is a constant, 
and the following {\em closed form equations 
of submanifolds $V_\beta$}: 
\begin{equation}\label{eq:3.58}
x = \int g(y)dy - C_1 y + C_2,  \;\; 
v = C_3  - \displaystyle \int \frac{du}{C_1 - f(u)}, 
\end{equation}
where $C_2$ and $C_3$ are constants. Hence we proved again 
that {\em the almost Grassmann structure with the base forms 
$(3.59)$ is $\beta$-semiintegrable.} In addition we proved that 
 {\em the family of submanifolds $V_\alpha$ on 
the manifold $M^4$ carrying the $AG (1, 3)$-structure 
with the basis forms $(3.50)$ depends 
on three parameters $C_1, C_2$ and $C_3$, and 
through any point $(x_0, y_0, u_0, v_0)$ 
of $M^4$ there passes a one-parameter family 
of submanifolds $V_\alpha$.}

If we look for conditions for $\alpha$-semiintegrability 
of this structure, then after writing equations (3.20) 
and taking their exterior derivatives we come to equations 
(3.39) which imply (3.40) and (3.41). So this structure is 
$\alpha$-semiintegrable and subsequently integrable if and only if 
 $f (u) = au + b$. Moreover,  the expression for $\lambda$ 
and the closed form 
equations of submanifolds $V_\alpha$ are (3.44) and (3.46). 

Note that if $g (y) = 0$, equations (3.57) and (3.58) 
coincide with equations (3.35) and (3.37).

Note that the application of formulas (3.18) gives the same values 
for 10 components of the tensor of conformal curvature 
that were obtained in Example {\bf 3.5}. 
}

Next we consider an example of an $\alpha$-semiintegrable 
almost Grassmann structure $AG (1, 3)$.

\examp{\label{examp:3.10} Suppose  that 
the basis 1-forms $\omega_\alpha^i$ of an almost Grassmann 
structure $AG (1, 3)$ are
\begin{equation}\label{eq:3.59}
\renewcommand{\arraystretch}{1.3}
\left\{
\begin{array}{ll}
\omega_1^3 = dx + p (y) du, & \omega_2^3 = dy, \\ 
\omega_1^4 = du,& \omega_2^4 = dv. 
\end{array}
\right.
\renewcommand{\arraystretch}{1}
\end{equation}

The metric of the $CO (2, 2)$-structure which is equivalent to 
the $AG (1, 3)$-structure with basis forms (3.59) is 
\begin{equation}\label{eq:3.60}
g = dx dv - dy du + p(y) du dv.
\end{equation}

In this case we have 
\begin{equation}\label{eq:3.61}
\renewcommand{\arraystretch}{1.3}
\left\{
\begin{array}{ll}
\omega_1^2 =  - \frac{1}{2} p' (y) \omega_1^4, \;\;\;\;
\omega_2^1 = 0,  \\
\omega_4^3 = - \frac{1}{2} p'(y) \omega_2^3, 
\;\;\;\;\omega_3^4 =  0,\\
\omega_1^1 = - \omega_2^2 
= \omega_3^3 = - \omega_4^4  = - \frac{1}{4} p' (y) \omega_2^4; 
\end{array}
\right.
\renewcommand{\arraystretch}{1}
\end{equation}
\begin{equation}\label{eq:3.62}
B_4 = - \frac{1}{4} p''(y), \;\; E_4 = - \frac{1}{4} (p' (y))^2;
\end{equation}
\begin{equation}\label{eq:3.63}
b_{134}^{122} = b_{134}^{212} = b_{134}^{221} 
= - \frac{1}{4} p''(y); 
\end{equation}
\begin{equation}\label{eq:3.64}
\renewcommand{\arraystretch}{1.3}
\left\{
\begin{array}{ll}
\omega_3^1 = 0,  &\omega_4^1 =  0, \\
\omega_3^2 = -\frac{1}{4} p''(y) \omega_2^4,  & 
\omega_4^2 = - \frac{1}{4} p''(y) \omega_2^3 
- \frac{1}{4} (p' (y))^2 \omega_2^4;
\end{array}
\right.
\renewcommand{\arraystretch}{1}
\end{equation}
\begin{equation}\label{eq:3.65}
\renewcommand{\arraystretch}{1.3}
\left\{
\begin{array}{ll}
\Omega_1^1 = \Omega_2^2 = = \frac{1}{4} p''(y) 
\omega_2^3 \wedge \omega_2^4, \; \;
\Omega_2^1  = \Omega_3^4 = \Omega_4^3 = \Omega_3^3 = \Omega_4^4 =0, \\
 \Omega_1^2 = - \frac{1}{4} p''(y) (\omega_1^3 \wedge \omega_2^4 
+  \omega_2^3 \wedge \omega_1^4).
\end{array}
\right.
\renewcommand{\arraystretch}{1}
\end{equation}

In this case  $b^{i\gamma\delta}_{jkl} = 0$, that is, 
$b_\alpha^1 = 0$, and {\em the almost Grassmann structure 
$AG (1, 3)$ with the basis forms $(3.59)$ 
is $\alpha$-semiintegrable}. However, in this case 
 $b_{\alpha kl}^{\beta \gamma \delta} \neq 0$ since 
 we have $b_{134}^{212} = - \frac{1}{4} p'' (y)$, that is
 $b_\alpha^1 \neq 0$, 
and in general {\em the almost Grassmann structure $AG (1, 3)$ 
with the basis forms $(3.59)$  is not $\beta$-semiintegrable, and 
consequently, it is not locally flat.}

If we apply the second method to the structure (3.59), 
we find that 
\begin{equation}\label{eq:3.66}
\lambda = - \frac{1}{p (y) + C_1},
\end{equation}
and the following {\em closed form 
equations of submanifolds $V_\alpha$}:
\begin{equation}\label{eq:3.67}
x = C_2-C_1 u,  \;\; v = C_3 + \displaystyle 
\int \frac{dy}{p (y) + C_3},
\end{equation}
where $C_1, C_2$ and $C_3$ are  constants.

Thus {\em the family of submanifolds $V_\alpha$ on 
the manifold $M^4$ carrying the $AG (1, 3)$-structure 
with the basis forms $(3.59)$ depends 
on three parameters $C_1, C_2$ and $C_3$, and 
through any point $(x_0, y_0, u_0, v_0)$ of $M^4$ there passes an 
one-parameter family of submanifolds $V_\alpha$.}

If  $p (y) = ay + b$, then we find that 
\begin{equation}\label{eq:3.68}
\mu = au + C_4, 
\end{equation}
and  the following {\em closed form 
equations of submanifolds $V_\beta$}: 
\begin{equation}\label{eq:3.69}
 x = C_5 - ayu - C_4 y - b u , \;\; 
 v = C_6 - \frac{1}{a} (au+C_4). 
\end{equation}
Thus {\em if $p (y) = ay + b$, then 
the family of submanifolds $V_\alpha$ on 
the manifold $M^4$ carrying the $AG (1, 3)$-structure 
with the basis forms $(3.49)$ depends 
on three parameters $C_4, C_5$ and $C_6$, and 
through any point $(x_0, y_0, u_0, v_0)$ of $M^4$ there passes a 
one-parameter family of submanifolds $V_\alpha$.}
}

Let us indicate examples of three more almost Grassmann 
structures similar to the structure (3.59). 

\examp{\label{examp:3.11} 
\begin{equation}\label{eq:3.70}
\renewcommand{\arraystretch}{1.3}
\left\{
\begin{array}{ll}
\omega_1^3 = dx, &\omega_2^3 = dy + q (x) dv, \\ 
 \omega_1^4 = du,& \omega_2^4 = dv. 
\end{array}
\right.
\renewcommand{\arraystretch}{1}
\end{equation}
}
\examp{\label{examp:3.12} 
\begin{equation}\label{eq:3.71}
\renewcommand{\arraystretch}{1.3}
\left\{
\begin{array}{ll}
\omega_1^3 = dx, & \omega_2^3 = dy,\\ 
\omega_1^4 = du + r (v) dx, & \omega_2^4 = dv. 
\end{array}
\right.
\renewcommand{\arraystretch}{1}
\end{equation}
}
\examp{\label{examp:3.13} 
\begin{equation}\label{eq:3.72}
\renewcommand{\arraystretch}{1.3}
\left\{
\begin{array}{ll}
\omega_1^3 = dx, &\omega_2^3 = dy, \\ 
\omega_1^4 = du, & \omega_2^4 = dv + s (u) dy. 
\end{array}
\right.
\renewcommand{\arraystretch}{1}
\end{equation}

As was the case for  the structure (3.59), {\em each of 
the structures $(3.70)-(3.72)$ is $\alpha$-semiintegrable and in 
general  not $\beta$-semiintegrable. These structures will be 
$\beta$-semiintegrable and consequently integrable 
if and only if the functions $q (x), r (v)$ and $s (u)$ are 
linear functions of $x, v$, and $u$, respectively.}
}

If we apply formulas (3.18) to Examples {\bf 3.10--3.13}, 
 we obtain that 
$C_\alpha = 0$ for all these examples (that is, the corresponding 
structures are $\alpha$-semiintegrable), and that the only 
nonvanishing components of the tensor $C_\beta$ are 
$b_3$ for Examples {\bf 3.10} ($b_3 = -\frac{p'' (y)}{4}$) 
and {\bf 3.12} ($b_3 =  \frac{r'' (v)}{4}$) 
and $b_1$ for Examples {\bf 3.11} ($b_1 = \frac{q'' (x)}{4}$) 
and {\bf 3.13} ($b_1 = -\frac{s'' (u)}{4}$). It follows that 
in the first case the polynomial 
$C_\beta$ has the triple root $\mu_1= \mu_2= \mu_3 
= \infty$ and the  simple root $\mu_4 = 0$, and in 
the second case 
it has the simple root $\mu_1 = \infty$  and the triple root 
$\mu_2 = \mu_3 = \mu_4 = 0$. 
According to Section {\bf 5.4} of the book \cite{AG96}, 
this means that the fiber bundle $E_\beta$ has two invariant 
distributions $\Delta_\beta (\infty)$ and $\Delta_\beta (0)$, and 
the distribution corresponding to a multiple root is integrable. 
Moreover, it is easy to prove that the distribution corresponding 
to a simple root is also integrable. 
For all cases the integral submanifolds $V_\beta$ are defined by 
the equations $y = C_2, v = C_4$ (for $\mu = \infty$) 
and by the equations $x = C_1, u = C_3$ (for $\mu = 0$).

{\bf 4.} 
An almost Grassmann structure is associated 
with a web, and if a web is transversally geodesic or 
isoclinic, the corresponding almost Grassmann structure 
is $\alpha$- or $\beta$-semiintegrable (see 
\cite{AG96}, Ch. 7). 
Our next two examples are generated by  examples of  
exceptional (nonalgebraizable) isoclinic webs of maximum 2-rank 
(see \cite{G85}, \cite{G86} or  \cite{G88}, Ch. 8, and 
\cite{AG96}, \S 5.5). 

\examp{\label{examp:3.14} Suppose  that 
the basis 1-forms $\omega_\alpha^i$ of an almost Grassmann 
structure $AG (1, 3)$ are
\begin{equation}\label{eq:3.73}
\renewcommand{\arraystretch}{1.3}
\left\{
\begin{array}{ll}
\omega_1^3 = dx, & \omega_2^3 = du,\\ 
\omega_1^4 = -(y+v)dx + (u-x)dy, & \omega_2^4 = (y+v) du + (u-x)dv. 
\end{array}
\right.
\renewcommand{\arraystretch}{1}
\end{equation}

The metric of the $CO (2, 2)$-structure which is equivalent to 
the $AG (1, 3)$-structure with basis forms (3.73) is 
\begin{equation}\label{eq:3.74}
g = 2(y + v) dx du + (u - x) (dx dv - dy du).
\end{equation}
In this case we have 
\begin{equation}\label{eq:3.75}
\renewcommand{\arraystretch}{1.3}
\left\{
\begin{array}{ll}
\omega_1^2 = \displaystyle \frac{1}{2(u-x)} \omega_1^3, &
\omega_3^4 = -  \displaystyle 
\frac{2(y+v)}{u-x} (\omega_1^3 + \omega_2^3) 
+  \displaystyle 
\frac{1}{2(u-x)} (-\omega_1^4 + \omega_2^4),   \\
\omega_2^1 = - \displaystyle \frac{1}{2(u-x)}  \omega_2^3, &
\omega_4^3 = 0,\\
\omega_1^1 = - \omega_2^2 =  \displaystyle 
\frac{3}{4(u-x)} (\omega_1^3 + \omega_2^3), &  
\omega_3^3 = - \omega_4^4  =  \displaystyle 
\frac{1}{4(u-x)} (-\omega_1^3 + \omega_2^3); 
\end{array}
\right.
\renewcommand{\arraystretch}{1}
\end{equation}
\begin{equation}\label{eq:3.76}
A_1 = B_2 = - \frac{5}{4(u-x)^2};\;\; 
A_2 = - \frac{7}{4(u-x)^2}. 
\end{equation}
\begin{equation}\label{eq:3.77}
b_{333}^{412} = - \displaystyle \frac{8(y+v)}{(u-x)^2}; 
\end{equation}
\begin{equation}\label{eq:3.78}
\renewcommand{\arraystretch}{1.3}
\left\{
\begin{array}{ll}
\omega_3^1 = -  \displaystyle 
\frac{1}{4(u-x)^2}(5\omega_1^3 + 7\omega_2^3),  &
\omega_4^1 = 0, \\
\omega_3^2 = -  \displaystyle 
\frac{1}{4(u-x)^2} (7\omega_1^3 + 5\omega_2^3),  & 
\omega_4^2 = 0;
\end{array}
\right.
\renewcommand{\arraystretch}{1}
\end{equation}
\begin{equation}\label{eq:3.79}
\renewcommand{\arraystretch}{1.3}
\left\{
\begin{array}{ll}
\Omega_1^2 = \Omega_2^1 =  \Omega_1^1 = \Omega_2^2 = 
\Omega_4^3 = \Omega_3^3 = \Omega_4^4 = 0, \\
\Omega_3^4 =  -  \displaystyle 
\frac{8(y+v)}{(u-x)^2} \omega_1^3 \wedge \omega_2^3. 
\end{array}
\right.
\renewcommand{\arraystretch}{1}
\end{equation}

It is easy to check that 
 $b_{\alpha kl}^{(\beta \gamma \delta)} = 0$, that is, 
$b_\beta^2 = 0$, and {\em the almost Grassmann structure 
$AG (1, 3)$ with basis forms $(3.73)$ is $\beta$-semiintegrable.} 
However, in this case 
 $b^{i \gamma \delta}_{(jkl)} \neq 0$. Thus, {\em the 
almost Grassmann structure $AG (1, 3)$ with basis forms $(3.73)$ 
is not locally flat.}

If we apply the second method for the structure (3.73), 
then while looking for $\alpha$-semiintegrability conditions, 
we find that 
$$
d\lambda = dy + dv - \frac{y+v+\lambda}{u-x} dx 
 + \frac{y+v-\lambda}{u-x} du.
$$
Exterior differentiation of this equation gives 
$\lambda dx \wedge du = 0$. Since 
on $V_\alpha$ we have $dx \wedge du \neq 0$, {\em the 
structure $(3.73)$ is never $\alpha$-semiintegrable.} 

For $\beta$-semiintegrability we find that 
\begin{equation}\label{eq:3.80}
\renewcommand{\arraystretch}{1.3}
\left\{
\begin{array}{ll}
d\mu = \displaystyle \frac{2 \mu (1 - \mu) du}{u -x}, \\
dx + \mu du = 0, \\
dy = - \mu dv - \displaystyle \frac{2 \mu (y + v)}{u-x} du.
\end{array}
\right.
\renewcommand{\arraystretch}{1}
\end{equation}
It is easy to see from (3.80) that 
$$
\frac{d \mu}{2 \mu (1 - \mu)} = \frac{du}{u - x} = 
\frac{dx}{- \mu (u - x)},
$$
and subsequently
\begin{equation}\label{eq:3.81}
\frac{d(u - x)}{(u - x)(1 + \mu)} =  \frac{d \mu}{2 \mu 
(1 - \mu)}.
\end{equation}
The solution of  equation (3.81)  is 
\begin{equation}\label{eq:3.82}
\frac{\sqrt{\mu}}{\mu-1} = C_1 (u-x).
\end{equation}
The submanifolds $V_\beta$ of the distribution 
$\Delta_\beta$ are defined by the  completely integrable 
system (3.80) where $\mu$ can be found from equations (3.82).

Thus {\em the family of submanifolds $V_\beta$ on 
the manifold $M^4$ carrying the $AG (1, 3)$-structure 
with the basis forms $(3.63)$ depends 
on three parameters, and 
through any point $(x_0, y_0, u_0, v_0)$ of $M^4$ there passes an 
one-parameter family of submanifolds $V_\alpha$.}

If we apply formulas (3.18), we find that for  
 the tensor of conformal curvature of the equivalent 
$CO (2, 2)$-structure we have $C_\beta = 0$ and the only 
nonvanishing component of $C_\alpha$ is $a_0 = - \frac{8(y+v)}{(u-x)^2}$. 
This means that the polynomial $C_\alpha (\lambda)$ has 
a quadruple root $\lambda = 0$, and the fiber bundle 
$E_\alpha$ possesses an integrable distribution $\Delta_\alpha (0)$. 
It is easy to prove that the submanifolds $V_\alpha$ of this
distribution are defined by the equations $x = C_1, u = C_3$, where $C_1$ and 
$C_3$ are constants.
}

\examp{\label{examp:3.15} Suppose  that 
the basis 1-forms $\omega_\alpha^i$ of an almost Grassmann 
structure $AG (1, 3)$ are
\begin{equation}\label{eq:3.83}
\renewcommand{\arraystretch}{1.3}
\left\{
\begin{array}{ll}
\omega_1^3 = dx, & \omega_2^3 = du,\\ 
\omega_1^4 = -v dx + u dy, & \omega_2^4 = y du - x dv. 
\end{array}
\right.
\renewcommand{\arraystretch}{1}
\end{equation}

The metric of the $CO (2, 2)$-structure which is equivalent to 
the $AG (1, 3)$-structure with basis forms (3.83) is 
\begin{equation}\label{eq:3.84}
g = (y + v) dx du - x dx dv - u dy du.
\end{equation}

In this case we have 
\begin{equation}\label{eq:3.85}
\renewcommand{\arraystretch}{1.3}
\left\{
\begin{array}{ll}
\omega_1^2 = - \displaystyle \frac{1}{2x} \omega_1^3, &
\omega_3^4 = \Bigl(\displaystyle \frac{y}{x} 
-  \displaystyle \frac{v}{u}\Bigr) 
(\omega_1^3 + \omega_2^3) -  \displaystyle 
\frac{1}{2u} \omega_1^4 -  \displaystyle 
\frac{1}{2x} \omega_2^4,   \\ 
\omega_2^1 = - \displaystyle \frac{1}{2u}  \omega_2^3, &
\omega_1^1 = - \omega_2^2 = \Bigl(\displaystyle \frac{1}{4u} 
-  \displaystyle \frac{1}{2x}\Bigr) 
\omega_1^3 +  \Bigl(\displaystyle \frac{1}{2u} 
-  \displaystyle \frac{1}{4x}\Bigr)\omega_2^3,\\
\omega_4^3 = 0,  &  
\omega_3^3 = - \omega_4^4  = \Bigl(\displaystyle 
\frac{1}{4u} +  \displaystyle \frac{1}{2x}\Bigr) 
\omega_1^3 +  \Bigl(\displaystyle \frac{1}{2u} 
+  \displaystyle \frac{1}{4x}\Bigr) \omega_2^3; 
\end{array}
\right.
\renewcommand{\arraystretch}{1}
\end{equation}
\begin{equation}\label{eq:3.86}
A_1 = \frac{1}{4u^2} + \frac{1}{4xu}; \;\;
A_2 = \frac{1}{4xu} - \frac{1}{4x^2} - \frac{1}{4u^2}; \;\;
B_2 = \frac{1}{2xu} - \frac{1}{4x^2}; 
\end{equation}
\begin{equation}\label{eq:3.87}
b_{333}^{412} = \Bigl(\displaystyle \frac{1}{u} 
- \displaystyle \frac{1}{x}\Bigr)
\Bigl(\displaystyle \frac{y}{x} 
- \displaystyle \frac{v}{u}\Bigr), \;\;
b_{333}^{312}= b_{334}^{412} =  - b_{343}^{412}
 =  \displaystyle \frac{1}{4x^2} 
-  \displaystyle \frac{1}{4u^2}; 
\end{equation}
\begin{equation}\label{eq:3.88}
\renewcommand{\arraystretch}{1.7}
\left\{
\begin{array}{ll}
\omega_3^1 = \Bigl(-\displaystyle \frac{1}{4u^2} + \frac{1}{2xu}\Bigr)
 \omega_1^3  
+   \displaystyle \Bigl( \frac{1}{4xu} - \frac{1}{4x^2} - \frac{1}{4u^2}
\Bigr) \omega_2^3,  &
\omega_4^1 = 0, \\
\omega_3^2 = \displaystyle \Bigl(\frac{1}{4xu} - \frac{1}{4x^2} 
- \frac{1}{4u^2} \Bigr) \omega_1^3 
+ \Bigl(\displaystyle \frac{1}{2xu}  - \displaystyle 
\frac{1}{4x^2}\Bigr) \omega_2^3,  & \omega_4^2 = 0;
\end{array}
\right.
\renewcommand{\arraystretch}{1}
\end{equation}
\begin{equation}\label{eq:3.89}
\renewcommand{\arraystretch}{1.7}
\left\{
\begin{array}{ll}
\Omega_1^1 =  \Omega_2^2 =  \Omega_1^2 = \Omega_2^1 = \Omega_4^3 = 0, \\
\Omega_3^3 = - \Omega_4^4 = \Bigl(\displaystyle \frac{1}{4u^2} 
-  \displaystyle \frac{1}{4x^2}\Bigr) 
\omega_1^3 \wedge \omega_2^3,\\
\Omega_3^4 = \Bigl(\displaystyle \frac{1}{u} 
- \displaystyle \frac{1}{x}\Bigr)\Bigl(\displaystyle \frac{y}{x} 
- \displaystyle \frac{v}{u}\Bigr) \omega_1^3 \wedge \omega_2^3 
+ \Bigl(\displaystyle \frac{1}{4x^2} - \displaystyle  \frac{1}{4u^2}\Bigr) 
(\omega_1^3 \wedge \omega_2^4  - \omega_2^3  \wedge \omega_1^4).
\end{array}
\right.
\renewcommand{\arraystretch}{1}
\end{equation}

It is easy to check that 
 $b_{\alpha kl}^{(\beta \gamma \delta)} = 0$, that is, 
$b_\beta^2 = 0$, and {\em the almost Grassmann structure 
$AG (1, 3)$ with the basis forms $(3.83)$  is $\beta$-semiintegrable.} 
However, in this case 
 $b^{i \gamma \delta}_{(jkl)} \neq 0$. Thus, {\em the 
almost Grassmann structure $AG (1, 3)$ with the basis forms 
$(3.83)$  is not locally flat.} 

If we apply the second method for the structure (3.73), 
then while looking for $\alpha$-semiintegrability conditions, 
we find that 
\begin{equation}\label{eq:3.90}
d\lambda = \frac{y+\lambda}{u} dx 
- dy + \frac{v-\lambda}{u} du + dv.
\end{equation}
Exterior differentiation of (3.90) gives 
\begin{equation}\label{eq:3.91}
\lambda = -\frac{1}{x^2+u^2}[vx(u-x) - yu(x+u)].
\end{equation}
The differential equations of submanifolds 
$V_\alpha$ have the form:
\begin{equation}\label{eq:3.92}
\renewcommand{\arraystretch}{1.3}
\left\{
\begin{array}{ll}
(v+y)(x+u) dx - (x^2+u^2) dy  = 0, \\
(v+y)(x+u) du - (x^2+u^2) dv  = 0.
\end{array}
\right.
\renewcommand{\arraystretch}{1}
\end{equation}
It is easy to prove that the integrability conditions of 
 equations (3.92) implies that  $du$ is proportional to $dx$. 
This is impossible since on $V_\alpha$ we have 
$dx \wedge du \neq 0$. Thus, {\em the structure 
$AG (1,3)$ with the basis forms $(3.83)$ is never 
$\alpha$-semiintegrable.} 

The $\beta$-submanifolds of this structure are 
defined by the following completely 
integrable system:
\begin{equation}\label{eq:3.93}
\renewcommand{\arraystretch}{1.3}
\left\{
\begin{array}{ll}
\mu du + dx = 0,  \\
\mu [(y+v) du - xdv] + udy  = 0,  \\
d \log (\mu - 1) 
= \Big(\displaystyle \frac{1}{u} 
- \displaystyle \frac{1}{x}\Bigr)dx.
\end{array}
\right.
\renewcommand{\arraystretch}{1}
\end{equation}

Thus {\em the family of submanifolds $V_\beta$ on 
the manifold $M^4$ carrying the $AG (1, 3)$-structure 
with the basis forms $(3.83)$ depends 
on three parameters, and 
through any point $(x_0, y_0, u_0, v_0)$ of $M^4$ there passes an 
one-parameter family of submanifolds $V_\beta$.}

If we apply formulas (3.18), we find that for  
 the tensor of conformal curvature of the equivalent 
$CO (2, 2)$-structure we have $C_\beta = 0$ and the only 
nonvanishing components of $C_\alpha$ are 
$a_0 = \Bigl(\frac{1}{u} - \frac{1}{x}\Bigr)\Bigl(\frac{y}{x} 
- \frac{v}{u}\Bigr)$ and $a_1 = \frac{1}{4}\Bigl(\frac{1}{x^2} 
- \frac{1}{u^2}\Bigr)$.   
This means that the polynomial $C_\alpha (\lambda)$ has 
a triple root $\lambda_1 = \lambda_2 = \lambda_3 = 0$ 
and a simple root $\lambda_4 = \frac{x+u}{xv - yu}$, and the fiber bundle 
$E_\alpha$ possesses an integrable distribution $\Delta_\alpha (0)$. 
It is easy to prove that the integral  submanifolds $V_\alpha$ of this
distribution are defined by the equations $x = C_1, u = C_3$, where 
$C_1$ and $C_3$ are constants. It is also easy to check that 
the distribution $\Delta (\lambda_4)$ is not integrable.
}

The next example is generated by an example of an isoclinic 
three-web given in  \cite{AS92} (see Exercise 6 on p. 133).

\examp{\label{examp:3.16} Suppose  that 
the basis 1-forms $\omega_\alpha^i$ of an almost Grassmann 
structure $AG (1, 3)$ are
\begin{equation}\label{eq:3.94}
\renewcommand{\arraystretch}{1.3}
\left\{
\begin{array}{ll}
\omega_1^3 = (y-v) dx + (x+u) dy, & \omega_2^3 = (y-v) du 
         - (x+u) dv, \\ 
\omega_1^4 = dy, & \omega_2^4 = dv. 
\end{array}
\right.
\renewcommand{\arraystretch}{1}
\end{equation}

The metric of the $CO (2, 2)$-structure which is equivalent to 
the $AG (1, 3)$-structure with basis forms (3.94) is 
\begin{equation}\label{eq:3.95}
g = (y - v) (dx dv - dy du) + 2 (x + u) dy dv.
\end{equation}

In this case we have 
\begin{equation}\label{eq:3.96}
\renewcommand{\arraystretch}{1.3}
\left\{
\begin{array}{ll}
\omega_1^2 = -\displaystyle \frac{1}{2(y-v)} \omega_1^4, &
\omega_3^4 = 0,   \\ 
\omega_2^1 = \displaystyle \frac{1}{2(y-v)}  \omega_2^3, &
\omega_4^3 = \displaystyle \frac{1}{2(y-v)} 
(\omega_1^3 - \omega_2^3) 
- \displaystyle \frac{2(x+u)}{y-v} (\omega_1^4 + \omega_2^4),\\
\omega_1^1 = - \omega_2^2 = -\displaystyle \frac{3}{4(y-v)} 
(\omega_1^4 +  \omega_2^4), &  
\omega_3^3 = - \omega_4^4  = -\displaystyle \frac{1}{4(y-v)} 
(\omega_1^4 -  \omega_2^4); 
\end{array}
\right.
\renewcommand{\arraystretch}{1}
\end{equation}
\begin{equation}\label{eq:3.97}
C_3 =- \displaystyle \frac{5}{4(y-v)^2}; \;\;
C_4 = -\displaystyle \frac{7}{4(y-v)^2}; \;\
E_4 = -\displaystyle \frac{5}{4(y-v)^2}; 
\end{equation}
\begin{equation}\label{eq:3.98}
b_{444}^{312} = \displaystyle \frac{8(x+u)}{(y-v)^2};
\end{equation}
\begin{equation}\label{eq:3.99}
\renewcommand{\arraystretch}{1.3}
\left\{
\begin{array}{ll}
\omega_3^1 = 0,  & \omega_4^1 = \displaystyle -\frac{1}{4(y-v)^2} 
(5\omega_1^4 + 7 \omega_2^4), \\
\omega_3^2 = 0,  & \omega_4^2 = \displaystyle -\frac{1}{4(y-v)^2} 
(7\omega_1^4 +5 \omega_2^4);
\end{array}
\right.
\renewcommand{\arraystretch}{1}
\end{equation}
\begin{equation}\label{eq:3.100}
\renewcommand{\arraystretch}{1.7}
\left\{
\begin{array}{ll}
\Omega_1^1 = \Omega_2^2 = \Omega_1^2 = \Omega_2^1 =  
\Omega_3^3 = \Omega_4^4 = \Omega_3^4 = 0,\\
\Omega_4^3 = \displaystyle \frac{8(x + u)}{(y-v)^2} 
\omega_1^4 \wedge \omega_2^4.
\end{array}
\right.
\renewcommand{\arraystretch}{1}
\end{equation}

It is easy to check that 
 $b_{\alpha kl}^{(\beta \gamma \delta)} = 0$, that is, 
$b_\beta = 0$, and the almost Grassmann structure 
$AG (1, 3)$ is $\beta$-semiintegrable. 
However, in this case 
 $b^{i \gamma \delta}_{(jkl)} \neq 0$. Thus, the 
almost Grassmann structure is not locally flat.

If we apply the method of direct integration, we find 
that the submanifolds $V_\beta$ are defined by the following 
completely integrable system of differential equations:
\begin{equation}\label{eq:3.101}
\renewcommand{\arraystretch}{1.3}
\left\{
\begin{array}{ll}
\displaystyle  \frac{d\mu}{2 \mu(\mu - 1)} 
= \displaystyle  \frac{dv}{y-v}, \\
(y-v) (\mu du + dx) - 2 \mu (x+u)dv = 0, \\
\mu dv + dy = 0.
\end{array}
\right.
\renewcommand{\arraystretch}{1}
\end{equation}

It is easy to see from (3.101) that 
$$
\frac{d \mu}{2 \mu (\mu - 1)} = \frac{dv}{y-v} = 
\frac{dy}{- \mu (y-v)},
$$
and subsequently
\begin{equation}\label{eq:3.102}
\frac{d(y - v)}{(y - v)(1 + \mu)} =  \frac{d \mu}{2 \mu 
(\mu - 1)}.
\end{equation}
The solution of  equation (3.102)  is 
\begin{equation}\label{eq:3.103}
\frac{\mu - 1}{\sqrt{\mu}} = C_1 (y-v).
\end{equation}

The submanifolds $V_\beta$ of the distribution 
$\Delta_\beta$ are defined by the  completely integrable 
system (3.101) where $\mu$ can be found from equations (3.103).

Thus {\em the family of submanifolds $V_\beta$ on 
the manifold $M^4$ carrying the $AG (1, 3)$-structure 
with the basis forms $(3.94)$ depends 
on three parameters, and 
through any point $(x_0, y_0, u_0, v_0)$ of $M^4$ there passes an 
one-parameter family of submanifolds $V_\beta$.}

If we apply formulas (3.18), we find that for  
 the tensor of conformal curvature of the equivalent 
$CO (2, 2)$-structure we have $C_\beta = 0$ and the only 
nonvanishing component of $C_\alpha$ is 
$a_4 = - \frac{8(x+u)}{(y-v)^2}$. 
This means that the polynomial $C_\alpha (\lambda)$ has 
a quadruple root $\lambda = \infty$, and the fiber bundle 
$E_\alpha$ possesses an integrable distribution $\Delta_\alpha (\infty)$. 
It is easy to prove that the submanifolds $V_\alpha$ of this
distribution are defined by the equations $y = C_2, v = C_4$, where 
$C_2$ and $C_4$ are constants.
}

{\bf 5.} The next example is generated by 
an example of an almost Grassmann 
structure $AG (1, 4)$ associated with a six-dimensional 
Bol web considered in \cite{AS92} (p. 270).

\examp{\label{examp:3.17} Suppose  that 
the basis 1-forms $\omega_\alpha^i$ of an almost Grassmann 
structure $AG (1, 4)$ are
\begin{equation}\label{eq:3.104}
\renewcommand{\arraystretch}{1.3}
\left\{
\begin{array}{ll}
\omega_1^3 = dx, & \omega_2^3 = du, \\ 
\omega_1^4 = [-2v - 2uz + (2x-1)w] e^{-2x}] dx + dy 
+ ue^{-2x} dz, & \omega_2^4 = (z du + dv - x dw)e^{-2x},\\ 
\omega_1^5 = dz, & \omega_2^5 = dw. 
\end{array}
\right.
\renewcommand{\arraystretch}{1}
\end{equation}

We prove that {\em the $AG(1, 4)$-structure  $(3.114)$ 
is $\alpha$-semiintegrable.} We apply the method of 
the direct integration.

For the structure (3.104) equations (3.22) 
of surfaces $V_\alpha, \dim V_\alpha = 2$ take the form 
\begin{equation}\label{eq:3.105}
\renewcommand{\arraystretch}{1.3}
\left\{
\begin{array}{ll}
\{\lambda_1 +  [-2v - 2uz + (2x-1)w] e^{-2x}\} dx 
+ dy + ue^{-2x} dz = 0, \\
\lambda_1 du + (z du + dv - x dw)e^{-2x} = 0, \\
 \lambda_2 dx +  dz = 0, \\
 \lambda_2 du +  dw = 0.
\end{array}
\right.
\renewcommand{\arraystretch}{1}
\end{equation}

Note that on $V_\alpha$ we have $dx \wedge du \neq 0$. 
Exterior differentiation of the last two equations 
of (3.105) gives the following exterior quadratic equations: 
$$
d \lambda_2 \wedge dx = 0, \;\;
d \lambda_2 \wedge du = 0.
$$
Since $dx \wedge du \neq 0$, it follows that $d \lambda_2 = 0$ 
and 
\begin{equation}\label{eq:3.106}
\lambda_2 = C_2,
\end{equation}
 where $C_2$ is a constant. Now the last two equations 
of (3.105) can be written as 
\begin{equation}\label{eq:3.107}
dz = - C_2 dx = 0, \;\; dw = -C_2 du.
\end{equation}
Integration of (3.107) gives 
\begin{equation}\label{eq:3.108}
z = - C_2 x +C_3, \;\; w = -C_2 u + C_4,
\end{equation}
where $C_3$ and $C_4$ are constants.

By (3.106) and (3.107), the first two equations of (3.105) become
\begin{equation}\label{eq:3.109}
\renewcommand{\arraystretch}{1.3}
\begin{array}{ll}
dy=\{- \lambda_1 + [2v + 2uz + (1- 2 x)w + C_2 u] e^{-2x}\} dx,\\ 
dv = -(\lambda_1 e^{2x} + z + C_2 x)du. 
\end{array}
\renewcommand{\arraystretch}{1}
\end{equation}
Exterior differentiation of (3.109) leads to the following 
exterior quadratic equations:
$$
(d \lambda_1  + 2\lambda_1 du) \wedge dx = 0, \;\; 
(d \lambda_1  + 2\lambda_1 dx) \wedge du = 0. 
$$
These equations can be written as 
$$
[d \lambda_1  + 2\lambda_1 (dx+du)] \wedge dx = 0, \;\; 
[d \lambda_1  + 2\lambda_1 (dx+du)] \wedge du = 0. 
$$
Since $dx \wedge du \neq 0$, it follows 
that 
$$
d \lambda_1  + 2\lambda_1 (dx+du) = 0.
$$
 This implies that 
\begin{equation}\label{eq:3.110}
\lambda_1 = C_1 e^{-2(x+u)}.
\end{equation}

Substituting this value of $\lambda_1$ and $z$ from 
(3.108) into the second equation of (3.109), we find that 
$$
dv = -(C_1 e^{-2u} + C_3) du.
$$
This implies that 
\begin{equation}\label{eq:3.111}
v = \frac{1}{2} C_1 e^{-2u} - C_3 u + C_5,
\end{equation}
where $C_5$ is a constant. 

Substituting the values of $\lambda_1, z$ and $w$ from 
(3.111), (3.110) and (3.108) into the first 
equation of (3.109), we find that 
$$
dy = [2C_5 + C_4 (1-2x)] e^{-2x} dx.
$$
The solution of this equation is
\begin{equation}\label{eq:3.112}
y = (- C_5  + C_4 x) e^{-2x} + C_6,
\end{equation}
where $C_6$ is a constant. 

Thus two-dimensional surfaces $V_\alpha$ are defined by the 
closed form equations (3.108), (3.110), (3.111) and (3.112). 
This completes the proof that the $AG(1, 4)$-structure  (3.104) 
is $\alpha$-semiintegrable. Thus {\em the family of submanifolds $V_\alpha$ 
on the manifold $M^6$ carrying the 
$AG (1, 4)$-structure with the basis forms $(3.104)$ depends 
on six parameters, and 
through any point $(x_0, y_0, z_0, u_0, v_0, w_0)$ 
of $M^6$ there passes a 
two-parameter family of submanifolds $V_\alpha$.}

We will prove now that {\em this structure is not 
$\beta$-semiintegrable.} 
Equations (3.23) of surfaces $V_\beta, \dim V_\beta = 3$, 
can be written as follows:
\begin{equation}\label{eq:3.113}
\mu du +  dx = 0, \;\; \mu \omega_2^4 +  \omega_1^4 = 0, \;\;
\mu dw + dz = 0.
\end{equation}
where $\mu = - \frac{s_1}{s_2}$. Taking exterior derivatives of 
the first and  third equations of (3.113), we find that 
$d\mu \wedge du = 0$ and $d\mu \wedge dw = 0$. Since on surfaces 
$V_\beta$ we must have $du \wedge dv \wedge dw \neq 0$, it 
follows that $d \mu = 0$, and 
\begin{equation}\label{eq:3.114}
\mu = C,
\end{equation}
where $C$ is a constant. Substituting this value of $\mu$ 
into the first and third equations of (3.113), we find that 
\begin{equation}\label{eq:3.115}
dx = - C du, \;\; dz = -C dw,
\end{equation}
and 
\begin{equation}\label{eq:3.116}
x = - C u + C_1, \;\; z = -C w + C_2,
\end{equation}
Substituting $\mu, x$ and $z$ from (3.114) and (3.115) 
into the second equation of (3.113), we obtain
\begin{equation}\label{eq:3.117}
\frac{1}{C} e^{2(C_1 - Cu)} = - dv + [C_1 + u(1-C)] dw 
+ [(C + 2C_1 - 1) w -2v -2C_2 u - C_2]du.
\end{equation}
Taking exterior derivative of (3.117), we find that if 
the surfaces $V_\beta$ exist, we would have 
\begin{equation}\label{eq:3.118}
 du \wedge [(2C - 1)dv + 2(1-C)(Cu + C_1 - 2)] dw] = 0.
\end{equation}
This is impossible since (3.118) implies that 
$du$ is a linear combination of $dv$ and $dw$ and as a result, 
$du \wedge dv \wedge dw = 0$. 

Thus  the $AG (1, 4)$-structure  is not 
$\beta$-semiintegrable. 
}

Example 3.5 can be generalized in the following manner.

\examp{\label{examp:3.18} Suppose that $x_\alpha, y_\alpha, 
\alpha = 1, \ldots , p$,  
are coordinates in $M^{2p}$, and that 
the basis 1-forms $\omega_\alpha^i$ of an almost Grassmann 
structure $AG (p-1, p+1)$ are
\begin{equation}\label{eq:3.119}
\renewcommand{\arraystretch}{1.3}
\left\{
\begin{array}{lll}
\omega_1^{p+1} = dx_1 + f(y_1) dx_2, & \omega_2^{p+1} = dx_2, &
     \omega_s^{p+1} = dx_s,\\ 
\omega_1^{p+2} = dy_1, & \omega_2^{p+2} = dy_2,  &
     \omega_s^{p+2} = dy_s, 
\end{array}
\right.
\renewcommand{\arraystretch}{1}
\end{equation}
where $s = 3, \ldots, p$. 
If the structure (3.119) is $\alpha$-semiintegrable, 
then $dy_1 \wedge dy_2 \wedge \ldots \wedge dy_p \neq 0$, and 
 the rows of the matrix $(\omega_\alpha^i)$ are proportional: 
\begin{equation}\label{eq:3.120}
 dx_1 + f(y) dx_2 + \lambda dy_1 = 0, \;\; 
dx_2 + \lambda dy_2 = 0, \;\;  dx_s + \lambda dy_s = 0. 
\end{equation}
Exterior differentiation of (3.120) gives the following exterior 
quadratic equations:
\begin{equation}\label{eq:3.121}
f' (y_1) dy_1 \wedge dx_2 + d\lambda \wedge dy_1 = 0, \;\; 
 d\lambda \wedge  dy_2 = 0, \;\;  d\lambda \wedge  dy_3 = 0.
\end{equation}
The last two equations of (3.121) imply that $d\lambda = 0$. 
By (3.120), the first equation of (3.121) gives 
$$
\lambda f' (y_1) dy_1 \wedge dy_2 = 0.
$$
Since $dy_1 \wedge dy_2 \neq 0$, it follows that 
$f' (y_1) = 0$ and $f (y_1) = a$, where 
$a$ is a constant. Thus the structure (3.119) 
 is $\alpha$-semiintegrable if and only if $f (y_1)$ is constant. 
If it is the case, closed form equations of integral submanifolds 
$V_\alpha$ of this structure have the form
\begin{equation}\label{eq:3.122}
x_1 + a x_2 + by_1 = C_1, \;\; 
x_2 + b y_2 = C_2, \;\;  x_s + b y_s = C_s.
\end{equation}
Hence   {\em the submanifolds $V_\alpha$ are flat 
$p$-dimensional submanifolds, and the family of submanifolds $V_\alpha$ 
depends on $p + 1$ constants $b, C_1, C_2$ and $C_s$.}

If the structure (3.129) is $\beta$-semiintegrable, 
then $dx_2 \wedge dy_2  \neq 0$, and 
 the columns of the matrix $(\omega_\alpha^i)$ are proportional: 
\begin{equation}\label{eq:3.123}
\renewcommand{\arraystretch}{1.3}
\left\{
\begin{array}{ll}
 dx_1 + f(y_1) dx_2 + \mu_1 dy_x = 0, &  dx_s + \mu_s dx_2 = 0, 
\\
dy_1 + \mu_1 dy_2 = 0, &  dy_s + \mu_s dy_2 = 0. 
\end{array}
\right.
\renewcommand{\arraystretch}{1}
\end{equation}
Exterior differentiation of (3.123) give the following exterior 
quadratic equations:
\begin{equation}\label{eq:3.124}
\renewcommand{\arraystretch}{1.3}
\left\{
\begin{array}{ll}
 (d\mu_1 + f' (y_1)) dy_1) \wedge dx_2 = 0, & 
d \mu_s \wedge  dx_2 = 0, \\  
(d\mu_1 + f' (y_1)) dy_1) \wedge  dy_2 = 0, & d\mu_s \wedge  
dy_2 = 0.
\end{array}
\right.
\renewcommand{\arraystretch}{1}
\end{equation}
Since $dx_2 \wedge dy_2 \neq 0$, it follows from 
(3.114) that 
$$
d\mu_1 + f' (y_1) dy_1 = 0, \;\; d\mu_s = 0,
$$
 and 
\begin{equation}\label{eq:3.125}
\mu_1 + f (y_1) = C_1, \;\; 
\mu_s = C_s.
\end{equation}
Substituting (3.125) into system (3.123), we find that 
\begin{equation}\label{eq:3.126}
\renewcommand{\arraystretch}{1.3}
\left\{
\begin{array}{ll}
 dx_1 + C_1 dx_2  = 0, &  dx_s + C_s dx_2 = 0, \\
dy_1 + (C_1 - f(y_1)) dy_2 = 0, &  dy_s + C_s dy_2 = 0. 
\end{array}
\right.
\renewcommand{\arraystretch}{1}
\end{equation}
The solution of this system is 
\begin{equation}\label{eq:3.127}
\renewcommand{\arraystretch}{1.3}
\left\{
\begin{array}{ll}
 x_1 + C_1 x_2  = A_1, &  x_s + C_s x_2 = A_s, \\
 \displaystyle \int \frac{dy_1}{C_1 - f(y_1)} + y_2 = B_1, & 
y_s + C_s y_2 = B_s. 
\end{array}
\right.
\renewcommand{\arraystretch}{1}
\end{equation}
Hence {\em the two-dimensional integral submanifolds $V_\beta$ 
are defined by closed  
form equations $(3.127)$, and the family of submanifolds 
$V_\beta$ depends on $3(p - 1)$ constants $C_1, A_1, B_1, C_s, 
A_s$ and $B_s$}.

Thus {\it if $f (y_1)$ is not a constant, the structure 
$AG (p-1, p+1)$ with the structure forms $(3.119)$ 
is $\beta$-semiintegrable, and if $f (y_1)$ is  a constant, this 
structure is locally flat}.
}

Example 3.10  can be also generalized.

\examp{\label{examp:3.19} Suppose that $x^3, x^4, x^t, y^3, 
y^4, y^t, \;t =  5, \ldots , q + 2$,  
are coordinates in $M^{2q}$, and that 
the basis 1-forms $\omega_\alpha^i$ of an almost Grassmann 
structure $AG (1, q+1)$ are
\begin{equation}\label{eq:3.128}
\renewcommand{\arraystretch}{1.3}
\left\{
\begin{array}{ll}
\omega_1^3 = dx^3 + p(y^3) dx^4, & \omega_2^{p+1} = dy^3,      \\ 
\omega_1^4 = dx^4, & \omega_2^4 = dy^4,  \\
 \omega_1^t = dx^t, &  \omega_2^t = dy^t,
\end{array}
\right.
\renewcommand{\arraystretch}{1}
\end{equation}
where $t = 5, \ldots, q+2$. 

If the structure (3.128) is $\alpha$-semiintegrable, 
then $dx^4 \wedge dy^4 \neq 0$ and 
 the rows of the matrix $(\omega_\alpha^i)$ are proportional: 
\begin{equation}\label{eq:3.129}
\renewcommand{\arraystretch}{1.3}
\left\{
\begin{array}{ll}
dx^3 + p (y^3) dx^4 + \lambda^3 dx^4 = 0, &  
dy^3 + \lambda^3 dy^4 = 0, \\  
dx^t + \lambda^t dx^4 = 0, & dy^t + \lambda^t dy^t = 0. 
\end{array}
\right.
\renewcommand{\arraystretch}{1}
\end{equation}
Exterior differentiation of (3.129) give the following exterior 
quadratic equations:
\begin{equation}\label{eq:3.130}
\renewcommand{\arraystretch}{1.3}
\left\{
\begin{array}{ll}
(d\lambda^3 + p' (y^3) dy^3) \wedge dx^4 = 0, &
 (d\lambda^3 + p' (y^3) dy^3)  \wedge  dy^4 = 0, \\
d\lambda^t \wedge dx^4 = 0, & d\lambda^t \wedge  dy^4 = 0. 
\end{array}
\right.
\renewcommand{\arraystretch}{1}
\end{equation}
Since $dy^3 \wedge dy^4 \neq 0$, equations  (3.130) imply that 
$d\lambda^3 + p' (y^3) dy^3=0$ and $d\lambda^t = 0$. 
Thus we have 
\begin{equation}\label{eq:3.131}
\lambda^3 + p (y^3) = C^3, \;\; \lambda^t = C^t.
\end{equation}
As a result, equations (3.129) become
\begin{equation}\label{eq:3.132}
\renewcommand{\arraystretch}{1.3}
\left\{
\begin{array}{ll}
dx^3 + C^3 dx^4  = 0, &  
dy^3 + (C^3 - p (y^3)) dy^4 = 0, \\  
dx^t + C^t dx^4 = 0, & dy^t + C^t dy^4 = 0. 
\end{array}
\right.
\renewcommand{\arraystretch}{1}
\end{equation}
It follows from (3.132) that closed form equations 
of two-dimensional  integral submanifolds $V_\alpha$ are 
\begin{equation}\label{eq:3.133}
\renewcommand{\arraystretch}{1.3}
\left\{
\begin{array}{ll}
 x^3 + C^3 x^4 = A^3, & \displaystyle   \int 
\frac{dy^3}{C^3 - p(y^3)} + y^4 
= B^4,\\
x^t + C^t x^4 = A^t, &  y^t + C^t y^4 = B^t. 
\end{array}
\right.
\renewcommand{\arraystretch}{1}
\end{equation}
Thus {\it the family of 
submanifolds $V_\alpha$ depends on $3 (q-1)$ 
constants $C^3, A^3, B^3, C^t, A^t$ and $B^t$.}

If the structure (3.128) is $\beta$-semiintegrable, 
then $dy^3 \wedge dy^4 \wedge \ldots \wedge dy^{q+2} \neq 0$, and 
 the columns of the matrix $(\omega_\alpha^i)$ are proportional: 
\begin{equation}\label{eq:3.134}
dx^3 + p (y^3) dx^4 + \mu dy^3 = 0, \;\; 
dx^4 +  \mu dy^4 = 0, \;\; 
dx^t +  \mu dy^t = 0.  
\end{equation}
Exterior differentiation of (3.134) gives the following exterior 
quadratic equations:
\begin{equation}\label{eq:3.135}
p' (y^3) dy^3 \wedge dx^4 + d\mu \wedge dy^3 = 0, \;\; 
 d\mu \wedge dy^4 = 0, \;\; d \mu \wedge  dy^t = 0.
\end{equation}
It follows from (3.135) that $d\mu = 0$ and
consequently $\mu = C_0$. 
As a result, by (3.134), the first equation of 
(3.135) becomes 
$$
C_0 p' (y^3) dy^3 \wedge dy^4 = 0.
$$
Since  $dy^3 \wedge dy^4 \neq 0$, it follows that 
$p' (y^3) = 0$ and $p (y^3) = a$. Thus {\em the structure $(3.118)$ is 
$\beta$-semiintegrable if and only if $p (y^3)$ is a 
constant.} If it is the case, equations (3.134) become 
\begin{equation}\label{eq:3.136}
dx^3 + a dx^4 + C_0 dy^3 = 0, \;\;  
dx^4 + C_0 dy^4 = 0, \;\;  
dx^t + C_0 dy^t = 0, 
\end{equation}
and closed form equations of integral submanifolds 
$V_\beta$ of this structure have the form
\begin{equation}\label{eq:3.137}
x^3 + a x^4 + C_0 y^3 = C^3, \;\;  
x^4 + C_0 y^4 = C^4, \;\;  
x^t + C_0 y^t = C^t. 
\end{equation}
Hence {\em the submanifolds $V_\beta$ are flat 
two-dimensional submanifolds, and the family of submanifolds 
$V_\beta$ depends on $q + 2$ constants $a, C^0, 
C^3, C^4$ and $C^t$.}

Thus {\it if $p (y^3)$ is not a constant, the structure 
$AG (1, q+1)$ with the structure forms $(3.128)$ 
is $\alpha$-semiintegrable, and if $p (y^3)$ is  a constant, this 
structure is locally flat}.
}

Example 3.17 can be generalized to an 
example of an $\alpha$-integrable 
almost Grassmann structure $AG (1, q +1)$. 

\examp{\label{examp:3.20} Suppose  that 
the basis 1-forms $\omega_\alpha^i$ of an almost Grassmann 
structure $AG (1, q+1)$ are
\begin{equation}\label{eq:3.138}
\renewcommand{\arraystretch}{1.3}
\left\{
\begin{array}{ll}
\omega_1^3 = dx, & \omega_2^3 = du,\\
\omega_1^4 = [-2v - 2uz + (2x-1)w] e^{-2x}] dx 
+ dy + ue^{-2x} dz, & \omega_2^4 = (z du + dv - x dw)e^{-2x},\\
\omega_1^s = dz^s, &  \omega_2^s = dw^s,  
\end{array}
\right.
\renewcommand{\arraystretch}{1}
\end{equation}
where $s = 5, \ldots, q + 2$.

Equations (3.105) will take the form

\begin{equation}\label{eq:3.139}
\renewcommand{\arraystretch}{1.3}
\left\{
\begin{array}{ll}
\{\lambda_1 +  [-2v - 2uz + (2x-1)w] e^{-2x}\} dx 
+ dy + ue^{-2x} dz = 0, \\
\lambda_1 du + (z du + dv - x dw)e^{-2x} = 0, \\
 \lambda_{s-3} dx +  dz^s = 0, \\
 \lambda_{s-3} du +  dw^s = 0, \;\; s = 5, \ldots, q + 2.
\end{array}
\right.
\renewcommand{\arraystretch}{1}
\end{equation}

The proof of the fact that  {\em the $AG(1, q+1)$-structure  
with the structure forms $(3.138)$ is $\alpha$-semiintegrable} is 
similar to that for Example 3.17.
}

Note that using this method of generalization, we cannot 
find a semiintegrable $AG (2, 5)$-structure 
(in this case $p = q = 3$). In fact, if we set 
$$
\renewcommand{\arraystretch}{1.3}
\left\{
\begin{array}{lll}
\omega_1^4 = \omega, & \omega_2^4 = dx_2, &\omega_3^4 = dx_3,\\    
\omega_1^5 = dy_1, & \omega_2^5 = dy_2, &\omega_3^5 = dy_3,\\    
\omega_1^6 = dz_1, & \omega_2^6 = dz_2, &\omega_3^6 = dz_3,    
\end{array}
\right.
\renewcommand{\arraystretch}{1}
$$
then for $\alpha$-semiintegrability we must have
$$
\renewcommand{\arraystretch}{1.3}
\left\{
\begin{array}{lll}
 \omega + \lambda dy_1 = 0, & dx_2 + \lambda dy_2 = 0, 
    &dx_3  + \lambda dy_3 = 0,\\    
 dz_1 + \widetilde{\lambda} dy_1 = 0, & 
 dz_2 + \widetilde{\lambda} dy_2 = 0,  &
 dz_3 + \widetilde{\lambda} dy_3 = 0.     
\end{array}
\right.
\renewcommand{\arraystretch}{1}
$$
It is easy to find from this that 
$\lambda = C, \;\; \widetilde{\lambda} 
= \widetilde{C}$, where $C$ and  $\widetilde{C}$ 
are constants, $ d\omega = 0$, that is, $\omega$ 
is a total differential, $\omega = dx_1$. But in this case 
the structure in question is also $\beta$-semiintegrable, 
i.e., this structure is locally flat.

\end{document}